\documentclass[a4paper]{article}

\usepackage[english]{babel}
\usepackage[utf8x]{inputenc}
\usepackage{amsmath}
\usepackage{graphicx}
\usepackage[colorinlistoftodos]{todonotes}

\usepackage{amsmath}
\usepackage{epsfig,amsthm}
\usepackage{latexsym}
\usepackage{amsfonts}
\usepackage{amssymb}
\usepackage{amscd}
\usepackage{mathrsfs}
\usepackage[all,cmtip]{xy}
\usepackage{enumerate}
\usepackage{tikz}

\usepackage{bbm}

\usepackage{color}
\usepackage[colorlinks]{hyperref}
\hypersetup{pdfstartview={FitH},         linkcolor=blue,  citecolor=green }

\numberwithin{equation}{section} 

\usepackage[a4paper,bindingoffset=0.2in,left=1in,right=1in,top=1in,bottom=1in,footskip=.25in]{geometry}

\newtheorem{thm}{Theorem}[section]
\newtheorem{cor}[thm]{Corollary}
\newtheorem{lem}[thm]{Lemma}
\newtheorem{prop}[thm]{Proposition}

\theoremstyle{definition}
\newtheorem{dfn}[thm]{Definition}
\newtheorem{rmk}[thm]{Remark}

\newcommand{\Fo}{{F_{\bullet}}}
\newcommand{\Eo}{{E_{\bullet}}}
\newcommand{\qo}{{q_{\bullet}}}

\newcommand{\tr}{{\mathrm{tr}}}

\newcommand{\GL}{{\mathrm{GL}}}
\newcommand{\Mat}{{\mathrm{Mat}}}

\newcommand{\SO}{{\mathrm{SO}}}

\newcommand{\sgn}{{\mathrm{sgn}}}

\newcommand{\Ind}{{\mathrm{Ind}}}
\newcommand{\cInd}{{\mathrm{cInd}}}
\newcommand{\Res}{{\mathrm{Res}}}

\newcommand{\diag}{{\mathrm{diag}}}

\newcommand{\Aut}{{\mathrm{Aut}}}

\newcommand{\Hom}{{\mathrm{Hom}}}
\newcommand{\End}{{\mathrm{End}}}

\newcommand{\Gal}{{\mathrm{Gal}}}

\title{Base change for ramified unitary groups: \\
the strongly ramified case}
\author{Corinne Blondel \and Geo Kam-Fai Tam}

\AtEndDocument{\bigskip
  \textsc{CNRS-IMJ-PRG, Universit{\'e} Paris Diderot, Case 7012, 75205 Paris Cedex 13, France.} \par  
  \textit{E-mail address:}  \texttt{corinne.blondel@imj-prg.fr} \par
  \addvspace{\medskipamount}
  \textsc{Department of Mathematics, IMAPP, Radboud University, Postbus 9010, 6500 GL Nijmegen, Netherlands.} \par 
    \textit{E-mail address:} \texttt{geotam@science.ru.nl} 
 }

\begin{document}
\maketitle

\begin{abstract}
We describe a special case of base change of certain supercuspidal representations from a {\it ramified} unitary group to a general linear group, both defined over a p-adic field of odd residual characteristic. Roughly speaking, we require the underlying stratum of a given supercuspidal representation to be skew maximal simple, and the field datum of this stratum to be of maximal degree, tamely ramified over the base field, and quadratic ramified over its subfield fixed by the Galois involution that defines the unitary group. The base change of this supercuspidal representation is described by a canonical lifting of its underlying simple character, together with the base change of the level-zero component of its inducing cuspidal type, modified by a sign attached to a quadratic Gauss sum defined by the internal structure of the simple character. To obtain this result, we study the reducibility points of a parabolic induction and the corresponding module over the affine Hecke algebra, defined by the covering type over the product of types of the given supercuspidal representation and of a candidate of its base change. \end{abstract}

\tableofcontents

\section{Introduction}
The local Langlands correspondence for a general linear group over a non-Archimedean local field $F$ is, roughly speaking, a parametrization of its irreducible admissible representations in terms of representations of the Weil-Deligne group of $F$. An extensive study of this correspondence involves certain important invariants, for example L- and epsilon-factors (see \cite{HT,Hen-simple} for the characteristic 0 case, \cite{LRS} for the positive characteristic case, and \cite{Scholze-LLC} for a recent proof). Using moreover the processes of automorphic induction and base change \cite{Hen-Herb,AC}, one obtains an explicit description of the correspondence in the essentially tame case \cite{BH-ET1,BH-ET2,BH-ET3}.

After \cite{Row-u3,Moeg-base-change,Arthur-new-book,Mok-unitary, unitary-inner}, we know that the local Langlands correspondence for classical groups is closely related to that for general linear groups, at least for discrete series representations: their L-packets are parametrized by multiplicity-free (conjugate-)self-dual representations of the Weil-Deligne group. Therefore, in principle, we can understand the correspondence for classical groups via general linear groups. For unitary groups, this theory is explained in the context of \emph{base change} in \cite{Moeg-base-change} (for special orthogonal groups and symplectic groups, the theory is sometimes called local transfer \cite{Shahidi-Cogdell-2014}). More precisely, the Langlands parameter of an L-packet of a discrete series of a unitary group is the same as that of its base change, a representation of a general linear group. Hence describing Langlands parameters for L-packets of discrete series is equivalent to describing their base changes. This approach was adopted in \cite{Adler-Lansky-unram,Adler-Lansky-ram,Blasco-u2,Blasco-u3} for supercuspidal representations of small unitary groups. 

In our present paper, we describe a special case of base change for ramified unitary groups, which  complements the previous method in \cite{tam-very-cusp} for describing the local Langlands correspondence for  packets of supercuspidal representations of unramified quasi-split unitary groups. In the previous case, we determined first the inertial class of the base change using the method developed in \cite{BHS} (for symplectic groups), and then the exact base change by using Asai L-functions \cite{shahidi-complementary,Goldberg-Siegel-case-unitary}. In the following paragraphs, we will first describe our new results and the methodology for ramified unitary groups, and then explain why the previous method fails.

Let $\Fo$ be a non-Archimedean local field of residual characteristic $p$ where $p$ is odd, $F/\Fo$ be a quadratic extension, $V$ be a vector space over $F$, and $\tilde{G}=\GL_F(V)$. Suppose that $V$ is equipped with a Hermitian form defined with respect to $F/\Fo$, and let $G$ be the corresponding 
unitary group, the fixed-point subgroup of the Hermitian involution $\sigma$ of $\tilde{G}$. Take a supercuspidal representation $\pi$ of $G$, compactly induced from a cuspidal type containing a {\em simple} character, or in other words, such that 
its underlying self-dual semi-simple stratum $\mathbf{s}$ is indeed simple (see \cite{BK,stevens-ss-char,Stevens-Miya} for the related definitions). Suppose that $E$ is the field datum of the stratum, which is invariant under the Galois involution defined by $\sigma$, with fixed-point subfield $\Eo$. In our paper, we require that
\begin{equation*}
[E:F]=\dim_F V
\quad \text{ and }\quad 
E/\Eo\text{ is quadratic ramified.}
\end{equation*}
We call this the \emph{strongly ramified} case. The latter condition actually forces $F/\Fo$ to be also quadratic ramified. For the final computation, we additionally assume that $E/F$ is \emph{tamely ramified}.

We refer the full detail of constructing supercuspidal representations by cuspidal types to \cite{BK,stevens-supercusp} and only specify that, under the above conditions, a cuspidal type is indeed a character (a representation of degree 1). Suppose that the cuspidal type of a supercuspidal representation $\pi$ of $G$ is of the form $\rho \kappa_0$, where $\kappa_0$ is the p-primary beta-extension of a simple character $\theta$, in the sense of \cite{BHS}, and $\rho$ is a  level zero character, which is a character of $\{\pm 1\}$ in our situation. We now take $\tilde{\theta}$ the self-dual simple character for $\tilde{G}$ whose  
restriction to the $\sigma$-fixed-point subgroup is the square of $\theta$ (see (\ref{simple-char-base-change})), and $\tilde{\boldsymbol{\kappa}}_0$ its unique self-dual p-primary beta-extension (see Proposition \ref{unique bold kappa}). We also take $\tilde{\boldsymbol{\rho}}$ a self-dual level zero character of $E^\times$ to form a cuspidal type $\tilde{\boldsymbol{\rho}}\tilde{\boldsymbol{\kappa}}_0$ and compactly induce it to a supercuspidal $\tilde{\pi}$ of $\tilde{G}$.

The following theorem is the main result of our present paper, giving the conditions for $\tilde{\pi}$ to be the base change of $\pi$, which is the only member in its L-packet in our case. Recall that we have assumed the `strongly ramified' condition and also that $E/F$ is tamely ramified.

\begin{thm}
\label{intro main theorem}
Under the above assumptions, suppose that the simple characters $\tilde{\theta}$ and $\theta$ are related as above (or see (\ref{simple-char-base-change})), and the level zero characters $\tilde{\rho}$ and $\rho$ are related by 
\begin{equation*}
\tilde{\rho}|_{\boldsymbol{\mu}_E}=\left(\frac{\cdot }{\boldsymbol{\mu}_E}\right)^{f(E/F)-1}\qquad\text{ and }\qquad 
\tilde{\rho}(\varpi_E)=\rho(-1)\epsilon_z^P(\varpi_E,\mathbf{s}),
\end{equation*}
where
\begin{itemize}

\item $\boldsymbol{\mu}_E$ is the subgroup of $E^\times$ of roots of unity with order coprime to $p$, and  $\left(\frac{\cdot }{\boldsymbol{\mu}_E}\right)$ is the quadratic character of 
$\boldsymbol{\mu}_E$;

\item $\varpi_E$ is a chosen uniformizer of $E$, and $\epsilon_z^P(\varpi_E,\mathbf{s})$ is a sign attached to a quadratic Gauss sum (see (\ref{the-quad-Gauss-sum})), defined by $\varpi_E$ and the simple stratum $\mathbf{s}$ associated to $\theta$.
\end{itemize}
Then $\tilde{\pi}$ is the base change of $\pi$.
\qed\end{thm}

For example, when $\dim V=1$, then $E=F$. The sign $\epsilon_z^P(\varpi_E,\mathbf{s})$ is 1, and the relation in Theorem \ref{intro main theorem} becomes 
$$\tilde{\rho}(x) = {\rho}(x{}^\sigma x)\qquad\text{for all }\qquad 
x\in F^\times,$$
 which is exactly the base change for characters of $\mathrm{U}_1$. The main idea here is that: when $\dim V>1$, we have to modify our base change formula for the level zero component by taking the internal structure of the simple character into account. 
This idea is also exhibited in the description of local Langlands correspondences for general linear groups in the essentially tame case \cite{BH-ET1,BH-ET2, BH-ET3}, for quasi-split unramified unitary groups \cite{tam-very-cusp}, and for symplectic groups \cite{BHS}. In all these cases, the modifications are incarnated by certain characters (called rectifiers in the first case and amending characters in the second case) of elliptic tori defined by the field data of simple characters.

We now describe the methodology for the theorem. The first step, following a series of papers of the first author \cite{Bl-Bl-SP4,Blondel-cover-propag,Blondel-Weil}, is to study the reducibility points of the parabolic induction
$$\tilde{\pi}|\det|^s\rtimes {\pi}:=\Ind_P^{G_W}(\tilde{\pi}|\det|^s\boxtimes \pi),\qquad s\in \mathbb{C},$$ 
where $W$ is the Hermitian space $W=V^-\oplus V\oplus V^+$, with each $V^-$ and $V^+$  just $V$ as a vector space and $V^-\oplus V^+$ hyperbolic, and so the unitary group $G_W$ contains a parabolic subgroup $P$ with Levi component $M=\tilde{G}\times G$. By \cite{silberger-special}, when this parabolic induction is reducible at a point $s\in \mathbb{R}_{\geq 0}$, then $\tilde{\pi}$ is conjugate-self-dual, in which case we can use \cite{Moeg-base-change} to show that indeed  
in our case, 
 $s=1$, and $\tilde{\pi}$ is the base change of $\pi$.

We can apply the theory of covering types \cite{BK-cover} to get a preliminary information about the reducibility: using the type $\lambda_P$ in $G_W$ covering $\tilde{\lambda}\boxtimes \lambda$ the product of a pair of types $\tilde{\lambda}$ and $\lambda$ for $\tilde{\pi}$ and $\pi$ respectively, together with the categorial equivalence
$$\mathcal{R}^{[M,\tilde{\pi}
\boxtimes \pi]}(G_W)\rightarrow \text{Mod-}\mathcal{H}(G_W,\lambda_P)$$
between the Bernstein component of the inertial class of $\tilde{\pi}
\boxtimes \pi$ in $G_W$ and the module category of the Hecke algebra $\mathcal{H}(G_W,\lambda_P)$, we can single out two candidates for $\tilde{\pi}$ of self-dual representations in the inertial class of the base change of $\pi$, the two differing from each other by an unramified character of a finite order, such that $\tilde{\pi}|\det|^s\rtimes {\pi}$ is reducible at a point in $\mathbb{R}_{\geq 0}$.

The second step is to further study the structure of the Hecke algebra $\mathcal{H}(G_W,\lambda_P)$ as well as its modules. By \cite{BK-cover,Stevens-Miya}, when $\tilde{\lambda}$ is conjugate-self-dual, this Hecke algebra has two generators, denoted by $T_y$ and $T_z$ in this paper, respectively satisfying a quadratic relation of the same form 
$$
 T_w^2 - (b_wT_w + c_w) = (T_w+\epsilon_w )(T_w-\epsilon_w q^{f(E/F)})=0,\,\qquad w\in \{y,z\},
$$
where each $\epsilon_w$ is a sign. By computing the eigenvalues of these two generators, or equivalently the coefficients $b_w$, on the module corresponding to $\tilde{\pi}|\det|^s\rtimes {\pi}$, we show that, when $\tilde{\pi}|\det|^s\rtimes {\pi}$ is reducible at $s=1$ and $\tilde{\rho}|_{\boldsymbol{\mu}_E}$ is given as in Theorem \ref{intro main theorem}, the signs $\epsilon_y$, $\epsilon_z$, and $\epsilon_y\epsilon_z$ are respectively $\rho(-1)$, $\epsilon_z^P(\varpi_E,\mathbf{s})$, and $\tilde{\rho}( \varpi_E)$ (modulo 
some irrelevant factors which cancel with each other in the comparison). The detailed version of this result can be derived from Theorem \ref{eigenvalues of Ty and Tz} and Corollary \ref{The two eigenvalues}.

It turns out that the coefficient  $b_y$ is easy to compute, while $b_z$ involves some calculations similar to \cite{Kim-Hecke-2, Kim-Hecke-3} for large $p$  
and  \cite{Bl-Bl-SP4} for $\mathrm{Sp}_4$, related to the structure of the simple stratum $\mathbf{s}$. For simplicity of our paper, we further assume that $E/F$ is tamely ramified. This is the condition assumed in \cite{Howe,MR-unitary} and the series \cite{BH-ET1,BH-ET2,BH-ET3} (see also the warning in \cite[(2.2.6)]{BK}), as well as in \cite{Yu-cusp} more generally. It also facilitates comparisons, by the second author in \cite{thesis, Tam-ETJLC}, between the essentially tame local Langlands correspondence for inner forms of general linear groups \cite{BH-ET1, BH-ET2, BH-ET3, BH-ETJLC} with the twisted endoscopy theory \cite{LS,KS}.

Finally, we show that $\tilde{\pi}$ is indeed the base change of $\pi$ using a finiteness result from M{\oe}glin for the possibilities of $\tilde{\pi}$ such that $\tilde{\pi}|\det|^s\rtimes {\pi}$ is reducible for some $s\in \mathbb{R}_{\geq 0}$ (\cite[4.Prop]{Moeg-base-change}, \cite[Th.3.2.1]{Moeg-endosc-L-param} for quasi-split groups, and M{\oe}glin-Renard \cite[8.3.5]{Moeglin-Renard-non-quasi-split} for non-quasi-split groups). The result is obtained from Arthur's endoscopic character relations \cite{Arthur-local-char-relation} and their generalizations in twisted endoscopy \cite[XI.]{ Moeglin-Waldspurger-Stabilisation-2}, which require that char($F$) is 0. It is possible that we could apply an approach analogous to \cite[Th. 2.5]{BHS} to compute the reducibilities of $\tilde{\pi}|\det|^s\rtimes {\pi}$ \emph{for all} $\tilde{\pi}$ and obtain a finiteness result similar to M{\oe}glin's without the characteristic requirement, but in our strongly ramified case we rather take the above shortcut using M{\oe}glin to simplify the discussion. (See \cite{Ganapathy-Varma} for a result on the local Langlands correspondence in positive characteristic
for split classical groups without using 
reducibilities of induced representations.)

We now briefly explain why the previous method in \cite{tam-very-cusp} fails in the strongly ramified case. According to \cite{Moeg-base-change}, there is a notion of parity of a conjugate-self-dual supercuspidal representation of $\GL_n$, either conjugate-orthogonal (+) or conjugate-symplectic ($-$) but not both. In the non-strongly ramified case, the two conjugate-self-dual candidates   have opposite parities because they differ by an 
unramified character $\tilde{\chi}$  such that $\tilde{\chi}\circ N_{E/F}$ is conjugate-symplectic. 
 In this case we can determine the correct base change between the two by computing their parities using Asai L-functions for example \cite{Hen-ext-sym}. However, in the strongly ramified case, $\tilde{\chi}\circ N_{E/F}$ is then conjugate-orthogonal, and so the two conjugate-self-dual candidates have the same parity (see the appendix in Section \ref{section appendix} for a detailed discussion). This explains why the previous method no longer works, and we have to rely on the complete structure of the modules over  
 the Hecke algebra.

\subsection{Acknowledgements}

To be added.

\subsection{Notations}

Let $\Fo$ be a non-Archimedean local field, with ring of integers $\mathfrak{o}_{\Fo}$, its maximal ideal $\mathfrak{p}_{\Fo}$, and residue field $\mathbf{k}_{\Fo}=\mathfrak{o}_{\Fo}/\mathfrak{p}_{\Fo}$ of cardinality $\qo$ and odd characteristic $p$. Let $F/\Fo$ be a quadratic extension, whose residue field $\mathbf{k}_{F}=\mathfrak{o}_{F}/\mathfrak{p}_{F}$ has $q$ elements, such that $q=\qo^2$ in the unramified case, and $q=\qo$ in the ramified case.
Let $\boldsymbol{\mu}_F$ be the subgroup of roots of unity of $F$ whose orders are coprime to $p$.

We denote $U_F:=\mathfrak{o}_F^\times$ and  $U^k_F:=1+\mathfrak{p}_F^\times$ for $k\in \mathbb{Z}_{\geq 1}$. If $N_{F/\Fo}:F\rightarrow \Fo$ is the norm homomorphism, we denote $U_{F/\Fo}:=\mathfrak{o}_F^\times\cap \ker N_{F/\Fo}$ and $U^k_{F/\Fo}:=U_F^k \cap \ker N_{F/\Fo}$ for $k\in \mathbb{Z}_{\geq 1}$. The Galois group $\Gal({F/\Fo})$ is generated by an involutive automorphism ${c}$.

If $\Lambda:\mathbb{Z}\rightarrow S$ is a sequence into a set $S$, we extend $\Lambda$ from $\mathbb{Z}$ to $\mathbb{R}$ by putting $\Lambda(r) = \Lambda(\lceil r \rceil)$ and $\Lambda(r_+) = \Lambda(\lceil r_+ \rceil)$, where $\lceil r \rceil$ and $\lceil r_+ \rceil$ are the smallest integers $\geq r$ and $>r$ respectively.

\section{Review of known results}

\subsection{Unitary groups}
\label{section unitary groups}

Let $V$ be an $F$-vector space. We denote $\tilde{A}=\tilde{A}_V=\End_F(V)$ and $\tilde{G}=\tilde{G}_V=\Aut_F(V)$. Suppose that $V$ is equipped with a non-degenerate $(F/\Fo,\epsilon)$-Hermitian  form $h=h_V$, where $\epsilon = \epsilon_V=\pm 1$. If $X\mapsto \bar{X} $ is the conjugate-adjoint on $\tilde{A}$ defined by $h$, we define the adjoint anti-involution ${}^\alpha X=-{}\bar{X}$. Note that 
\begin{equation*}
{}^\alpha(XY)=-{}^\alpha Y{}^\alpha X,\,
\qquad
X,Y\in \tilde{A}.
\end{equation*}
We also have the corresponding involution 
$\sigma: X\mapsto \bar{X}^{-1}$
on $\tilde{{G}}$. The subgroup $G=\tilde{G}^\sigma$ is a (connected) unitary group that we consider throughout the paper and whose Lie algebra is ${A}=\tilde{A}^\alpha$.

Given an $\mathfrak{o}_F$-lattice $L$ in $V$, we denote by $L^*$ its dual  $\{v\in V:h(v,L)\subset \mathfrak{p}_F\}$. 
We call an $\mathfrak{o}_F$-lattice sequence $\Lambda$ in $V$ self-dual if there exists $d\in \mathbb{Z}$ such that $\Lambda(k)^*=\Lambda(d-k)$  for all $k\in \mathbb{Z}$. As in \cite{stevens-supercusp}, we always normalize a self-dual lattice sequence such that 
\begin{equation}
\label{middle-1-and-period-even}
\text{$d=1$ and its $\mathfrak{o}_F$-period $e(\Lambda/\mathfrak{o}_F)$ is even.}
\end{equation}

Given an $\mathfrak{o}_F$-lattice sequence $\Lambda$ in $V$, we define an $\mathfrak{o}_F$-lattice sequence $\tilde{\mathfrak{P}}^k_\Lambda$, for $k\in \mathbb{Z}$, in $\tilde{A}$ by
$$\tilde{\mathfrak{P}}^k_\Lambda = \{x\in \tilde{A}: x\Lambda(m)\subseteq \Lambda(m+k)\text{ for all }m\in \mathbb{Z}\}.$$
Hence $\tilde{\mathfrak{A}}_\Lambda:=\tilde{\mathfrak{P}}^0_\Lambda$ is a hereditary order in $\tilde{A}$, with Jacobson radical 
$\tilde{\mathfrak{P}}_\Lambda:=\tilde{\mathfrak{P}}^1_\Lambda$. We denote by $v_\Lambda$ the valuation on $\tilde{A}$ associated to $\Lambda$. If $\Lambda$ is self-dual, then each $\tilde{\mathfrak{P}}^k_\Lambda$ is $\alpha$-invariant, in which case we put ${\mathfrak{P}}_\Lambda=\tilde{\mathfrak{P}}_\Lambda^\alpha$.

We also define $\tilde{U}_\Lambda=\tilde{U}^0_\Lambda=\tilde{\mathfrak{A}}_\Lambda^\times$ 
and $\tilde{U}^k_\Lambda=1+\tilde{\mathfrak{P}}^k_\Lambda $ for $k\in \mathbb{Z}_{\geq 1}$. If $\Lambda$ is self-dual, then each $\tilde{U}^k_\Lambda$ is $\sigma$-invariant, in which case we put $U_\Lambda=\tilde{U}_\Lambda^\sigma$ and $U^k_\Lambda=(\tilde{U}^k_\Lambda)^\sigma$. The quotient $\mathsf{G}_{\Lambda}=U_\Lambda/U^{1}_\Lambda$ is the group of $\mathbf{k}_{\Fo}$-points of a reductive group $\boldsymbol{\mathsf{G}}_\Lambda$ defined over $\mathbf{k}_{\Fo}$, which is disconnected in general. We denote by $U^0_\Lambda$ the inverse image of ${\mathsf{G}}^0_{\Lambda}$ in $U_\Lambda$, where ${\mathsf{G}}^0_{\Lambda}$ is the subgroup of $\mathbf{k}_{\Fo}$-points of the identity component $\boldsymbol{\mathsf{G}}^0_{\Lambda}$ of $\boldsymbol{\mathsf{G}}_{\Lambda}$. 

Suppose that $\Lambda$ is a self-dual lattice sequence of the form 
\begin{equation*}
    \Lambda(0) \supseteq \Lambda(1) = \Lambda(0)^* \supseteq \Lambda(2) = \mathfrak{p}_F \Lambda(0),
\end{equation*}
then $U_\Lambda$ is a maximal compact subgroup, and $U^0_\Lambda$ is the underlying maximal parahoric subgroup. When $F/\Fo$ is unramified, $U^0_\Lambda$ is a product of at most two  unitary groups relative to $\mathbf{k}_{F}/\mathbf{k}_{\Fo}$; while when $F/\Fo$ is ramified it is a product of a symplectic group and a special orthogonal group, both defined over $\mathbf{k}_{F}=\mathbf{k}_{\Fo}$ and can be possibly trivial. 

\subsection{Cuspidal types}
\label{section cuspidal types}

We recall from \cite{BK,stevens-supercusp} the constructions of cuspidal types for general linear groups and unitary groups. The compact inductions of these types are irreducible supercuspidal representations.

Let $\mathbf{s}=[\Lambda,r,0,\beta]$ be either a simple stratum or the null stratum $[\Lambda,0,0,0]$ where, in the former case, we denote $E=F[
\beta]$, which is a field contained in $\title{A}$, such that  $\Lambda$ is an $\mathfrak{o}_E$-lattice chain, while in the latter case, we put $E=F$ and $\Lambda(k) = \Mat_n(\mathfrak{p}_F^k)$. We denote by  $\tilde{A}_{E}$ and $\tilde{G}_{E}$ respectively the centralizer of $\beta$ in $\tilde{A}_{}$ and $\tilde{G}_{}$, and for $k\in \mathbb{Z}$, denote $\tilde{\mathfrak{P}}^k_{\Lambda,{E}}
=\tilde{\mathfrak{P}}^k_\Lambda\cap \tilde{A}_{E}$ 
and $\tilde{U}^k_{\Lambda,E}=
\tilde{U}^k_{\Lambda_{}}\cap \tilde{G}_{E}$. As in \cite[Chapter 3]{BK} associated to $\mathbf{s}$ we construct 
\begin{itemize}
\item subrings $\tilde{\mathfrak{H}}=\tilde{\mathfrak{H}}_{\Lambda,\beta}\subseteq \tilde{\mathfrak{J}}=\tilde{\mathfrak{J}}_{\Lambda,\beta} $ of $\tilde{A}$ and the two-sided fractional ideals $\tilde{\mathfrak{H}}^k  =\tilde{\mathfrak{H}}\cap \tilde{\mathfrak{P}}^k_\Lambda$ and $\tilde{\mathfrak{J}}^k =\tilde{\mathfrak{J}}\cap \tilde{\mathfrak{P}}^k_\Lambda$,  for all $k\in \mathbb{Z}$;

\item subgroups $
\tilde{H}^{k} = \tilde{H}^{k}_{\Lambda,\beta}=\tilde{\mathfrak{H}}\cap \tilde{U}^k_\Lambda
\subset 
\tilde{J}^{k} = \tilde{J}^{k}_{\Lambda,\beta}=\tilde{\mathfrak{J}}\cap \tilde{U}^k_\Lambda$ of  $\tilde{{G}}$, for all $k\in \mathbb{Z}_{>0}$;

\item $\tilde{\mathcal{C}}(\mathbf{s}) := \tilde{\mathcal{C}}(\Lambda,0,\beta)$ the set of simple characters of $\tilde{H}^1$;

\item associated to each simple character $\tilde{\theta}\in \tilde{\mathcal{C}}(\mathbf{s})$  the Heisenberg representation $\tilde{\eta}=\tilde{\eta}_{\tilde{\theta}}$ of $
\tilde{J}^{1}$, an irreducible representation that restricts to a multiple of $\tilde{\theta}$ on $
\tilde{H}^{1}$;

\item a beta-extension $\tilde{\kappa}$ of $\tilde{\eta}$ to the subgroup $\tilde{J}=\tilde{J}_{\Lambda, \beta}=\tilde{U}_{\Lambda,{E}}\tilde{J}^1$ (note that among these extensions there is a unique one $\tilde{\kappa}_0$ whose determinant has a $p$-power order, called the p-primary beta-extension);

\end{itemize}

To construct a supercuspidal representation of $\tilde{G}$, we require that $\Lambda$ is principal, which is assumed from now on, and also $e(\Lambda/\mathfrak{o}_E)=2$ (note the convention in (\ref{middle-1-and-period-even})). We take an irreducible representation $\tilde{\rho}$ of $\tilde{J}$ inflated from a cuspidal representation of $\tilde{J}/\tilde{J}^1 \cong \GL_{n/[E:F]}(\mathbf{k}_{E})$. The maximal simple type $\tilde{\lambda}=\tilde{\kappa}\otimes \tilde{\rho}$ extends to an irreducible representation $\tilde{\boldsymbol{\lambda}}$ of $\tilde{\mathbf{J}}=E^\times \tilde{J}$, which is then compactly induced to an irreducible supercuspidal representation $\tilde{\pi}=\cInd_{\tilde{\mathbf{J}}}^{\tilde{G}}\tilde{\boldsymbol{\lambda}}$.

Note that if $\mathbf{s}$ is null, then by convention we assume that $e(\Lambda/\mathfrak{o}_F)=2$, and $\tilde{\mathcal{C}}(\mathbf{s})$ contains only the trivial character of $\tilde{H^1}=\tilde{J}^1=\tilde{U}^1_\Lambda$. We take $\tilde{\kappa}$ of $\tilde{U}_\Lambda$ to be trivial, and choose $\tilde{\rho}$ to be inflated from a cuspidal representation of $\tilde{U}_\Lambda/\tilde{U}^1_\Lambda\cong \GL_n(\mathbf{k}_F)$. The extension $\tilde{\boldsymbol\lambda}$ on $\tilde{\mathbf{J}}=F^\times \tilde
J$ of $\tilde{\lambda}=\tilde{\rho}$ can be chosen by fixing a central character, and so   
$\tilde{\pi}=\cInd_{\tilde{\mathbf{J}}}^{\tilde{G}}\tilde{\boldsymbol{\lambda}}$ is a level zero supercuspidal representation.

Every supercuspidal representation of $\tilde{G}$ is obtained by the way above. Moreover, the maximal simple type $\tilde{\lambda}$ is determined by $\tilde{\pi}$ up to conjugacy, and so is the extended type $\tilde{\boldsymbol{\lambda}}$ containing $\tilde{\lambda}$.

If $V$ is a Hermitian space, as in Section \ref{section unitary groups}, we call a simple stratum $\mathbf{s}$ skew if $\Lambda$ is self-dual and $\beta\in A$. In this case, all subgroups $\tilde{H}^1$, $\tilde{J}^1$, $\tilde{J}$, and $\tilde{\mathbf{J}}$ are $\sigma$-invariant. If we choose $\tilde \theta$, $\tilde{\kappa}$, $\tilde{\rho}$, and the extension of $\tilde{\boldsymbol{\lambda}}$ to be self-dual, i.e., $\sigma$-invariant, then so is $\tilde{\pi}$. Note that the p-primary beta-extension $\tilde{\kappa}_0$ of $\tilde{\theta}$ is self-dual because of its uniqueness. 
 If $\beta\neq0$, then $-\alpha$ restricts to a Galois involution on $E$. We denote the fixed field by $E_\bullet$. Since $\tilde{A}_{E}$ and  ${\tilde{\mathfrak{P}}^k_{\Lambda,{E}}}$, resp. $\tilde{G}_{E}$, $\tilde{U}_{\Lambda,{E}}$, and $\tilde{U}^k_{\Lambda,{E}}$, are invariant under $\alpha$, resp.  $\sigma$, we define ${A}_{E/\Eo}=\tilde{A}_{E}^\alpha$, ${G}_{E/\Eo} = \tilde{G}_{E}^\sigma$, and ${\mathfrak{P}^k_{\Lambda,{E/\Eo}}}$, $U_{\Lambda,{E/\Eo}}$,  ${U^k_{\Lambda,{E/\Eo}}}$ similarly.

When $\mathbf{s}$ is skew and semisimple, following \cite{stevens-supercusp}, we construct

\begin{itemize}

\item 
 subgroups $H^1:=(\tilde{H}^1)^\sigma\subseteq J^1:=(\tilde{J}^1)^\sigma\subset  {J}:=\tilde{J}^\sigma$;

\item  a set ${\mathcal{C}}(\mathbf{s}):={\mathcal{C}}(\Lambda,0,\beta)$ of semisimple characters of ${H}^1$, defined as follows: denote the subset of $\tilde{\mathcal{C}}(\mathbf{s})$ of self-dual semisimple characters by $\tilde{\mathcal{C}}(\mathbf{s})^\sigma$, and define ${\mathcal{C}}(\mathbf{s})$ to be the image set of the following map (well-defined since the group $H^1$ is a pro-$p$ subgroup where $p$ is odd)
\begin{equation}
\label{simple-char-base-change}
\tilde{\mathcal{C}}(\mathbf{s})^\sigma \rightarrow {\mathcal{C}}(\mathbf{s}),\,\qquad \tilde{\theta}\mapsto \theta:=(\tilde{\theta}|_{H^1})^{1/2}, 
\end{equation}
which turns out to be bijective since the restriction operation satisfies the properties of Glaubermann correspondence (see \cite[Sec 5.3]{stevens-supercusp} for details); 

\item associated to each semisimple character ${\theta}\in {\mathcal{C}}(\mathbf{s})$  the Heisenberg representation ${\eta}={\eta}_{{\theta}}$ of $
{J}^{1}$, an irreducible representation that restricts to a multiple of ${\theta}$ on $
{H}^{1}$;

\item a beta-extension ${\kappa}$ of ${\eta}$ to the subgroup ${J}= {U}_{\Lambda,{E/\Eo}} {J}^1$
 (note that among these extensions there is a unique one ${\kappa}_0$ whose determinant has a $p$-power order, called the p-primary beta-extension).

\end{itemize}

To construct a supercuspidal representation of ${G}$, from now on we suppose that ${G}_{E/\Eo}$ has a compact center and the parahoric subgroup $U^0_{\Lambda,{E/\Eo}}$ is maximal in $G_{E/\Eo}$. We take an irreducible representation ${\rho}$ of ${J}$ inflated from a cuspidal representation ${\bar{\rho}}$ of $\mathsf{G} = {J}/{J}^1 $, a (possibly non-connected) classical group over the finite field $\mathbf{k}_{\Fo}$. (A cuspidal representation of a disconnected $\mathsf{G}$ means that its restriction to $\mathsf{G}^0$ is a sum of conjugates of a cuspidal representation.) The product $\lambda=\kappa \otimes \rho$ is then a cuspidal type, and is compactly induced to an 
 irreducible supercuspidal representation ${\pi}=\mathrm{cInd}_{J}^{{G}}\lambda$. If $\mathbf{s}$ is null,
   then by definition ${\mathcal{C}}(\mathbf{s})$ contains only the trivial character of $U^1_\Lambda$. We again take $\kappa$ to be trivial, and so  
${\pi}=\mathrm{cInd}_{{{J}}}^{{G}}{\rho}$ is a level zero supercuspidal representation.

\begin{rmk}
${G}_{E/\Eo}$ has a compact center except precisely when it has a factor isomorphic to 
the split $\SO_2$, (which is just a $\GL_1$), which does not happen in our case. Also, for the parahoric subgroup $U^0_{\Lambda,{E/\Eo}}$ to be maximal, it is not enough to just assume that the order $\mathfrak{P}^0_{\Lambda,{E}}$ is maximal, since there exists such a maximal order whose corresponding parahoric subgroup is not maximal. See \cite[Appendix]{Stevens-Miya} for details.
\qed\end{rmk}

\subsection{Covering types}
\label{section covering types}

To proceed, we require a larger unitary group $G_W$ defined on the space $W=V\perp Z$. Here $Z=Z_-\oplus Z_+$ is an $\epsilon$-Hermitian space, with $Z_+$ being a finite dimensional vector space over $F$, and $Z_-$ is the dual space of $Z_+$ with respect to the form
\begin{equation*}
h_Z((z-,z+),(w_-,w_+)) = (z_-,w_+) +\epsilon\cdot      {}^c(w_-,z_+),
\quad 
z_-,w_-\in Z_-\text{ and }z_+,w_+\in Z_+.
\end{equation*}
The $\epsilon$-Hermitian form on $W$ is $h_W=h\perp h_Z$. We denote the Lie algebra of $G_W$ by $A_W$, the subspace of fixed points of $\tilde{A}_W=\End_F(W)$ by the adjoint anti-involution defined by $h_W$.

We denote by $M $ the Levi subgroup of block diagonal matrices in $G_W$, isomorphic to $\tilde{G}_{Z_-}\times G_V$. We denote by $P$ the subgroup of block upper triangular matrices 
leaving invariant the flag $Z_- \subset V\subset W$, 
with unipotent radical $U$ and the opposite $U^-$. We denote by $i_M:M\rightarrow G_W$ the embedding 
$(g,h)\mapsto \diag(g,h, {}^\sigma g)$ 
for $g\in \tilde{G},\,h\in G$, and abbreviate a matrix of the form 
$\left[\begin{smallmatrix}
I&X&Y \\ &I&{}^\alpha X \\ &&I
\end{smallmatrix}\right]\in U$ by $(X,Y)^+$, and similarly $(X,Y)^-=\left[\begin{smallmatrix}
I&& \\ {}^\alpha X &I&\\ Y &X&I
\end{smallmatrix}\right]\in U^-$. We can check that these matrices satisfy the relation
\begin{equation}
\label{relation between X and Y}
X{}^\alpha X=Y-{}^\alpha Y.
\end{equation}

If $\Lambda_{+}$ is an $\mathfrak{o}_F$-lattice chain in $Z_+$, we define
\begin{equation*}
\Lambda_{-}(k) = \{z_-\in Z_-:(z_-,\Lambda_+(1-k))\subset\mathfrak{p}_F\},\,\quad 
k\in \mathbb{Z}.
\end{equation*}
 If $\beta_+\in \tilde{A}_{Z_+}$, we define $\beta_-\in \tilde{A}_{Z_-}$ by 
\begin{equation*}
(\beta_- z_-,z_+)=- ^c(z_-,\beta_+z_+),\,
\quad\text{for }z-\in Z_-,z+\in Z_+.
\end{equation*}

We will proceed to construct covering types for $G_W$. To simplify the discussions in this paper, as well as to focus on the case we are interested 
in, let's assume the following  from now on.
\begin{itemize}
\item $V=Z_+$, so that we can identify $Z_-$ with $Z_+$ by the Hermitian form $h=h_V$. 

\item $\beta_+ = \beta = \beta_-$, in particular it means that $\beta$ is self-dual: ${}^\alpha\beta = \beta$. 

\item $\Lambda = \Lambda_+$, which is an $\mathfrak{o}_E$-lattice chain, with $E=F[\beta]$. We moreover assume that $\Lambda_V$ is principal and its  $\mathfrak{o}_{E}$-period is 2  (hence its $\mathfrak{o}_F$-period is  $2e(E/F)$).

\end{itemize}
We define a self-dual $\mathfrak{o}_{E}$-lattice sequence $\mathfrak{m}$ in $W=V \perp Z$ of $\mathfrak{o}_{E}$-period 6 (hence of $\mathfrak{o}_F$-period $6e(E/F)$) by  
\begin{equation*}
\mathfrak{m}(k)=\Lambda_-\left(\frac{k-1}{3}\right)
\oplus 
\Lambda\left(\frac{k}{3}\right)
\oplus 
\Lambda_+\left(\frac{k+1}{3}\right), \qquad k\in \mathbb{Z},
\end{equation*}
and two maximal self-dual $\mathfrak{o}_{E}$-lattice sequences $\mathfrak{M}^{w}$ in $W$, where $w\in \{y,z\}$ and both of $\mathfrak{o}_{E}$-period $2$, such that the set of lattices in the sequence $ \mathfrak{M}^{y}    $ is \begin{equation*}
\{\Lambda_-(k)\oplus \Lambda(k)\oplus \Lambda_+(k) : k \in \mathbb Z\},
\end{equation*}
while those in $ \mathfrak{M}^{z}   $ is   \begin{equation*}
\{\Lambda_-(k)\oplus \Lambda(k)\oplus \Lambda_+(k+1) ,\, \Lambda_-(k)\oplus \Lambda(k+1)\oplus \Lambda_+(k+1) : k \in \mathbb Z\}.
\end{equation*}

Given a skew simple or null stratum $\mathbf{s}=[\Lambda,r,0,\beta]$, we hence have corresponding skew semi-simple strata  
$\mathbf{s}_\mathfrak{m}=[\mathfrak{m},r_\mathfrak{m},0,\beta]$ and $\mathbf{s}_{w}=[\mathfrak{M}^w,r_{w},0,\beta]$,  
with $w\in \{y,z\}$, where $\beta$ is  embedded into $ \tilde{A}_W$ as the block-diagonal matrix with diagonal blocks $(\beta,\beta,\beta)$. 
We now follow \cite{BK,stevens-supercusp} to
define, for $\mathcal{L}=\mathfrak{m}$ or $\mathfrak{M}^w$, with $w\in \{y,z\}$,

\begin{itemize}

\item subgroups $H^1_{\mathfrak{L}}\subseteq J^1_{\mathfrak{L}}\subset J_{\mathfrak{L}}$ in $G_W$;

\item semi-simple characters $\theta_\mathfrak{L}$ on $H^1_{\mathfrak{L}}$, such that they satisfy a transfer relation with each other, as in \cite[(3.6)]{BK}, \cite[Sec 3.5]{stevens-ss-char};

\item beta-extensions $\kappa_{\mathfrak{M}^w}$ of $\theta_{\mathfrak{M}^w}$ on $J_{\mathfrak{M}^w}$, and $\kappa_{\mathfrak{m}}^{w}$ of $\theta_{\mathfrak{m}}$ on $J_{\mathfrak{m}}$, such that $\kappa_{\mathfrak{M}^w}$ and $\kappa_{\mathfrak{m}}^{w}$ are compatible in the sense of \cite[Lemma 4.3]{stevens-supercusp};

\item unique p-primary beta extension $(J_{\mathfrak{M}^w}, \kappa_{\mathfrak{M}^w,0})$, and the one $(J_{\mathfrak{m}},\kappa_{\mathfrak{m},0}^{w})$ compatible with $\kappa_{\mathfrak{M}^w,0}$ (note that $\kappa_{\mathfrak{m},0}^{w}$ is in general \emph{not} the p-primary beta extension $ \kappa_{\mathfrak{m},0}$ of $\theta_{\mathfrak{m}}$).

\end{itemize}

Since the decomposition $W= V\perp(Z_-\oplus Z_+)$ is properly subordinate to the stratum $\mathbf{s}_\mathfrak{m}$, the subgroups $H^1_{\mathfrak{m}}$, 
 $J^1_{\mathfrak{m}}$, and $J_{\mathfrak{m}}$
  admit an Iwahori decomposition of the form
\begin{equation*}
K=K^-K_MK^+,\quad\text{where }\quad K_{M}=K\cap M,\,K^+=K\cap U\text{, and }K^-=K\cap U^-, 
\end{equation*}
for $K$ being one of the compact subgroups just mentioned, and the product can be taken in any order. We can then show from \cite[Cor 5.11 and Prop 5.5]{stevens-supercusp} that $H^1_{M}:=\tilde{H}^1_{}\times H^1_{}$ and
\begin{equation*}
\theta_{\mathfrak{m}}|_{H^1_{M}}=\tilde{\theta}_{\Lambda_+}\boxtimes\theta_\Lambda
\end{equation*}
where $\tilde{\theta}_{\Lambda_+}$ and $\theta_\Lambda$ are simple characters related under the correspondence (\ref{simple-char-base-change}), i.e., $\theta_\Lambda = (\tilde{\theta}_{\Lambda_+}|_{H^1})^{1/2}$. %

To construct a covering type, we define
\begin{itemize}
\item the following subgroups
 \begin{equation*}
\begin{split}
H^1_{P}&=\tilde{H}^1_{\mathfrak{m}}(\tilde{J}^1_{\mathfrak{m}}\cap U)\cap G_W
=
{H}^1_{\mathfrak{m}}({J}^1_{\mathfrak{m}}\cap U),
\\
J^1_{P}&=\tilde{H}^1_{\mathfrak{m}}(\tilde{J}^1_{\mathfrak{m}}\cap P)\cap G_W
=
{H}^1_{\mathfrak{m}}({J}^1_{\mathfrak{m}}\cap P)
\text{, and }
\\
J_{P}&=\tilde{H}^1_{\mathfrak{m}}(\tilde{J}_{\mathfrak{m}}\cap P)\cap G_W
=
{H}^1_{\mathfrak{m}}({J}_{\mathfrak{m}}\cap P),
\end{split}
\end{equation*}
and also denote $J_P^\pm  = J_P\cap U^\pm  $;

\item a simple character $\theta_P$ on $H^1_{P}$ by extending $\theta_{\mathfrak{m}}$ trivially to ${J}^1_{\mathfrak{m}}\cap U$, and the unique representation $\eta_P$ on $J^1_P$ containing $\theta_P$ (by \cite[Lemma 5.12]{stevens-supercusp});

\item for $w\in \{y,z\}$, 
an extension $\kappa_P^w$ of  $\eta_P$ to $J^1_P$ that is the restriction of  $\kappa_{\mathfrak{m}}^w$ to the subspace of  $({J}^1_{\mathfrak{m}}\cap U)$-fixed vectors of $\eta_{\mathfrak{m}}$ (which can be also characterized by the relation $\kappa_{\mathfrak{m}}^{w}\cong \Ind_{J_{P}}^{J_{\mathfrak{m}}}\kappa^w_P$).

\end{itemize} 

With our construction of $\mathfrak{m}$, the decomposition $W= V\perp (Z_- \oplus Z_+)$ is moreover exactly subordinate to $\mathbf{s}_\mathfrak{m}$. By \cite[Prop 6.3]{stevens-supercusp}, if we write $\kappa^w_P|{J_M} = \tilde{\kappa}^w\boxtimes \kappa^w $, then $\tilde{\kappa}^w$ and $\kappa^w$ are beta-extensions of $\tilde{\theta}$ and $\theta$ respectively, and $\tilde{\kappa}^w$ is self-dual (with respect to $\sigma$)  \cite[Cor 6.10]{stevens-supercusp}.

We then choose a (product of) cuspidal representation ${\rho}_M=\tilde{\rho}\boxtimes \rho$ of $J_P$, inflated from $J_P/J_P^1\cong J_M/J_M^1$. (At this moment, there is no relation between $\tilde{\rho}$ and $\rho$.) Define 
$$\lambda_P^w=  \kappa_P^w \otimes {\rho}_M$$
which is a covering type over $\lambda_M^w:=\lambda_P^w|_{J_M}$ by \cite{Stevens-Miya}. Finally, we denote by $\pi_M^w=\tilde{\pi}^w\boxtimes \pi^w$ a supercuspidal representation of $M=\tilde{G}_Z\times G_V$ containing the maximal type $\lambda_M^w=\tilde{\lambda}^w\boxtimes \lambda^w$, where $\tilde{\lambda}^w=\tilde{\kappa}^w\otimes \tilde{\rho}$ and $\lambda^w= {\kappa}^w\otimes {\rho}$.

\subsection{Covers and Hecke algebras}
\label{section Covers and Hecke algebras}

At the beginning of this subsection, we use $G$ to denote the $F$-points of a connected reductive group over $F$, and will later switch back to denote a unitary group as in previous subsections. We recall the following notions from \cite{BK-cover}.
\begin{itemize}
\item  Suppose that $G$ contains a parabolic subgroup $P$ with Levi component $M$. If $\pi_M$ is a supercuspidal representation of $M$,
we denote by $\mathfrak{s}=[M,\pi_M]_G$ 
the inertial class of $\pi_M$, and by $\mathcal{R}^{\mathfrak{s}}(G) $ the full subcategory of representations of $G$ whose irreducible subquotients are those of 
the normalized parabolic induction $\iota_P^{G}\tau: = \Ind_P^G(\tau\Delta_P^{1/2})$, where $\tau$ ranges over representations in $\mathfrak{s}$ and $\Delta_P$ is the modulus  character of $P$.

\item Suppose that $K$ is a compact open subgroup of $G$. We fix a Haar measure on $G$ such that $K$ has volume $1$.  
Given a representation $\lambda$ of $K$ on a finite dimensional $\mathbb{C}$-vector space $W$, we denote by $\mathcal{H}(G,\lambda)$ the associated Hecke algebra, which is the space of compactly supported functions $f:G\rightarrow\mathrm{End}_{\mathbb{C}}(W)$ satisfying 
$$f(k_1gk_2)=\lambda(k_1)\circ f(g)\circ \lambda(k_2)\text{, for all }k_1,\,k_2\in K\text{ and }g\in G,$$
with an associative $\mathbb{C}$-algebra structure under the convolution
$$f_1*f_2(g)=\int_Gf_1(x)f_2(x^{-1}g)dx\text{, for all }g\in G.$$
The support of every element in $\mathcal{H}(G,\lambda)$ lies in the intertwining set $I_G(\lambda):=\{g\in G:\Hom_{K\cap gKg^{-1}}(\lambda,{}^g\lambda)\neq 0\}$.

\end{itemize}
We call $(K,\lambda)$ an $\mathfrak{s}$-type if, for an irreducible representation $\tau$ of $G$, 
\begin{equation*}
\tau|_K\text{ contains }\lambda\quad\Leftrightarrow \quad \text{the inertial class of the cuspidal support of $\tau$ is $\mathfrak{s}$},
\end{equation*}
in which case we have an equivalence of categories
\begin{equation*}
\mathcal{M}:\mathcal{R}^{\mathfrak{s}}(G)\rightarrow \text{Mod-}\mathcal{H}(G,\lambda),\,\qquad\tau\mapsto \Hom_{K}(\lambda,\tau).
\end{equation*}

We now return to our notations in the previous subsection, suppose that $\pi_M$ contains a type $(J_M,\lambda_M)$, and $(J_P,\lambda_P)$ is a  $G_W$-cover of $(J_M,\lambda_M)$ (e.g. $\lambda_P=\lambda^w_P$, for $w\in \{y,z\}$). In particular, we have $\lambda_P|_{J_M} = \lambda_M$. Since $(J_M,\lambda_M)$ is an $[M,\pi_M]$-type, we obtain 

\begin{itemize}
\item that $(J_P,\lambda_P)$ is a $[M,\pi_M]_{G_W}$-type, i.e., we have an equivalence of categories
\begin{equation*}
\mathcal{M}:\mathcal{R}^{[M,\pi_M]}(G_W)\rightarrow \text{Mod-}\mathcal{H}(G_W,\lambda_P),\,\qquad\tau\mapsto \Hom_{J_P}(\lambda_P,\tau);
\end{equation*}

\item an injective morphism of algebras \cite[(8.3, 8.4)]{BK-cover} 
\begin{equation*}
t_P:\mathcal{H}(M,\lambda_M)\rightarrow
\mathcal{H}(G_W,\lambda_P)
\end{equation*}
giving rise to the natural functor 
\begin{equation*}
    \begin{split}
        (t_P)_*:\text{Mod-}\mathcal{H}(M,\lambda_M)&\rightarrow\text{Mod-}\mathcal{H}(G_W,\lambda_P),\, \qquad
        \\
        D&\mapsto \Hom_{\mathcal{H}(M,\lambda_M)}(\mathcal{H}(G_W,\lambda_P),D),
    \end{split}
\end{equation*}
such that the following diagram
\begin{equation}\label{commutative diagram}
  \xymatrixcolsep{5pc}\xymatrix{
\mathcal{R}^{[M,\pi_M]}(G_W)
\ar[r]^{\mathcal{M}_{G_W}}    
&\text{Mod-}\mathcal{H}(G_W,\lambda_P)
\\
\mathcal{R}^{[M,\pi_M]}(M) \ar[u]_{\iota_P^{G_W}}
\ar[r]^{\mathcal{M}_M }
&\text{Mod-}\mathcal{H}(M,\lambda_P|_{J_M})  
\ar[u]_{(t_P)_*}
}
\end{equation}
commutes.

\end{itemize}

We will be interested in when the parabolic induction 
\begin{equation}
\label{definition of parabolic induction}
\tilde{\pi}|\det|^s\rtimes \pi  := \iota_P^{G_W}(\tilde{\pi}|\det|^s\boxtimes \pi),\,\qquad
s\in \mathbb{C},
\end{equation}
is reducible. By \cite{silberger-special}, when this happens, $\tilde{\pi}$ must be self-dual, and there is a unique \emph{real} $s_{\tilde{\pi},{\pi}}\geq 0$ such that (\ref{definition of parabolic induction}) is reducible at $s=\pm s_{\tilde{\pi},{\pi}}$. A simple argument shows that there are exactly two self-dual representations in the inertial class of $\tilde{\pi}$ giving this reducibility. Say one of those is $\tilde{\pi}$ and the 
 other $\tilde{\pi}'$, a twist  
of $\tilde{\pi}$ by an unramified character,
and denote by $ f_{\tilde{\pi}} = f_{\tilde{\pi}'}$ the order of the stabilizer of $\tilde{\pi}$ in the group of unramified characters of $\tilde{G}$. The points of reducibility for $\tilde{\pi}$ are of the form \begin{equation}
\label{4 reducibility points}
\left\{\pm s_1,\,\pm s_2+\frac{\pi i}{f_{\tilde{\pi}}\log q}\right\}\qquad\text{ for some }s_1,\,s_2\in \mathbb{R}_{\geq 0},
\end{equation}
and those for $\tilde{\pi}'$ takes the same form, with $s_1$ and $s_2$ exchanged. To obtain more information about these values, we will study the structures of Hecke algebras and their modules in the next section.

\subsection{Structures of Hecke algebras}
\label{section Structures of Hecke algebras}

We continue from the constructions in Subsection \ref{section Covers and Hecke algebras}, and briefly recall the structure of the above Hecke algebras, referring the detail to \cite{BK-cover,stevens-supercusp,Stevens-Miya}.

For $\tilde{G}$, if $E=F[\beta]$ is the field datum of $\mathbf{s}$, then $I_{\tilde{G}}(\tilde{\lambda})=  \tilde{\mathbf{J}}$ with $
\tilde{\mathbf{J}}/\tilde{J} \cong E^\times/U_E \cong \left<\varpi_E\right>$. The Hecke algebra $\mathcal{H}(\tilde{G},\tilde{\lambda})$ is isomorphic to $\mathbb{C}[Z,Z^{-1}]$, where $Z$, chosen up to a scalar, is supported on a single coset $\varpi_E \tilde{J}_\Lambda$, where $\varpi_E $ is a uniformizer of $E$ (note that this coset is independent of the choice of $\varpi_E $). For ${G}$, since $I_{{G}}({\lambda})={J}_\Lambda$, we have $\mathcal{H}({G},{\lambda})\cong \mathbb{C}$. Therefore, $\mathcal{H}(M,\tilde{\lambda}\boxtimes \lambda)\cong \mathbb{C}[Z,Z^{-1}]$.

We now fix $w$ to be either one of $\{y,z\}$ and abbreviate $\kappa_P=\kappa^w_P$ s and $\lambda_P=\kappa_P\otimes \rho_M$, such that $\lambda_P|_{J_M}=\tilde{\lambda}\boxtimes \lambda$. Consider the Hecke algebra $\mathcal{H}(G_W,\lambda_P)$. If $\tilde{\lambda}$ is not self-dual, 
then 
$$\mathcal{H}(G_W,\lambda_P)\cong \mathcal{H}(M,\tilde{\lambda}\boxtimes \lambda)\cong \mathbb{C}[Z,Z^{-1}].$$
 The interesting case happens when $\tilde{\lambda}$ is self-dual, in which case it is known that   \cite[Cor 6.16]{stevens-supercusp} \cite[Prop 3.3]{Blondel-Weil}     
\begin{equation*}
\mathrm{rank}_{\mathcal{H}(M,\tilde{\lambda}\boxtimes \lambda)}
(\mathcal{H}(G_W,\lambda_P))=
\#(N_{G_W}([M,\tilde{\pi}\boxtimes\pi])/M)=2.
\end{equation*}
The Hecke algebra $\mathcal{H}(G_W,\lambda_P)$ has two invertible generators $T_w$, for $w\in \{y,z\}$. To describe them precisely, we choose two elements $s_y$ and $s_z$ (for example, we may choose $s_1$ and $s^\varpi_1$ in \cite{stevens-supercusp} or \cite{Blondel-Weil}), each of which is a generator for the normalizer group $N_{G_W}([M,\tilde{\pi}\boxtimes\pi])$ mod $M$ such that 
\begin{equation*}
s_yJ_P^-s_y\subset J_P^+\text{ and }s_zJ_P^+s_z^{-1}\subset J_P^-,
\end{equation*}
and moreover
$
\zeta:=s_ys_z=i_M(\varpi_EI,I)
$
which is a $P$-positive element in the sense that
       $$\zeta J_P^+\zeta^{-1} \subset J_P^+\text{ and }\zeta^{-1}J_P^-\zeta\subset J_P^-.$$
(If we have chosen $P^-$ instead of $P$ to define our covering type, then we have to switch $y$ and $z$.) Each of the two generators $T_w$, for $w\in \{y,z\}$,   has support on a single double coset $J_Ps_wJ_P$ and is defined up to a scalar. They satisfy  
\begin{equation}
\label{product Ty and Tz}
T_y*T_z=t_P(Z),
\end{equation}
and certain quadratic relations
\begin{equation}
\label{quadratic relation general b and c}
T_w*T_w=b_wT_w+c_w     \mathbbm{1}  
\end{equation}
for certain real numbers $b_w$ and $c_w$. Here $\mathbbm{1}$ is the unit in $\mathcal{H}(G_W,\lambda_P)$, which is the function supported  on $J_P$ with $\mathbbm{1}(1)=I_{\lambda_P}$, the identity operator on the representation space of $\lambda_P$.

For direct formulas for the coefficients, we follow \cite[Sec 1]{Bl-Bl-SP4} and obtain
\begin{equation}
\label{formula for b and c}
\begin{split}
c_y  &=(\dim \lambda_P)^{-1} [J_P^+:s_yJ_P^- s_y]      \tr (T_y(s_y)T_y(s_y^{-1})),\, 
\\
b_y&=(\dim \lambda_P)^{-1}\sum_{u\in \frac{s_yJ_P^+s_y\cap J_Ps_yJ_P}{J_P^- }} {\tr T_y(u)},
\end{split}
\end{equation}
and similarly
\begin{equation}\label{formula for b and c for y }
\begin{split}
    c_z  &=(\dim \lambda_P)^{-1} [s_z J_P^-s_z^{-1}:J_P^+]  \tr (T_z(s_z)T_z(s_z^{-1})),\,
\\
b_z&=(\dim \lambda_P)^{-1}\sum_{u\in \frac{s_z J_P^-s_z^{-1}\cap J_Ps_zJ_P}{J_P^+}} {\tr T_z(u)}{} .
\end{split}
\end{equation}
Later in (\ref{sy and sz}) we will choose $s_w$ such that $s_w^2\in J_M$, so that we may normalize each $T_w$, up to a sign,
 such that
\begin{equation}
\label{T(s^2)}
T_w(s_w)^2=\lambda_P(s_w^2)
\end{equation}
which is equivalent to requiring that
\begin{equation*}
T_w(s_w) T_w(s_w^{-1})=I_{\lambda_P}
\end{equation*}
 and so that $c_y=[J_P^+:s_yJ_P^- s_y]    $ and $c_z = [s_z J_P^-s_z^{-1}:J_P^+] $ are both positive.

Indeed, by \cite[Section 7.1]{stevens-supercusp}, each $T_w$ comes from the generator of $\mathcal{H}(U_{\mathfrak{M}^w,E},\tilde{\rho}\boxtimes \rho)$  which is then reduced to a Hecke algebra for the finite reductive quotient (not necessarily connected)
$$\mathsf{G}_w  = U_{\mathfrak{M}^w,E} / U^1_{\mathfrak{M}^w,E}.$$
By Lusztig, under suitable normalizations, the quadratic relations can be written as 
\begin{equation}
\label{normalized quadratic relation general b and c}
(T_w+\mathbbm{1})*(T_w -q_E^{r_w}\mathbbm{1})=0.
\end{equation}
for certain integers $r_w\geq 0$. The values of $r_w$ can be determined by the method in \cite{Lust-Stevens,BHS} combined with Lusztig's classification \cite{lusz-classical}. In our present situation, we can also compute these values by comparing the coefficients, or more precisely the eigenvalues, between (\ref{normalized quadratic relation general b and c}) and (\ref{quadratic relation general b and c}), which will be done in Section \ref{section First coefficient}.

\subsection{Modules over Hecke algebras}
\label{section Modules of Hecke algebras}

It is worthwhile to mention a result, in   (\ref{values of real parts}) below, from \cite{Blondel-Weil} about the real parts of the points of reducibility, i.e., the values of $s_1$ and $s_2$ in (\ref{4 reducibility points}), which are enough to determine the inertial Jordan blocks (c.f. \cite{BHS}) of a supercuspidal $\pi$ of $G$ with type $\lambda$. The result will be recalled    after some prerequisites on modules over Hecke algebras.

 Suppose that $\tilde{\pi}=\cInd_{\tilde{\mathbf{J}}}^{\tilde{G}}\tilde{\boldsymbol{\lambda}}$, where $\tilde{\boldsymbol{\lambda}}$ extends $\tilde{\lambda}$ of $\tilde{J}$. We abbreviate $\mathcal{H}_{\tilde{G}}=\mathcal{H}(\tilde{G},\tilde{\lambda})$ and recall the equivalence of categories
$$\mathcal{M}_{\tilde{G}}:\mathcal{R}^{[\tilde{G},\tilde{\pi}]}(\tilde{G})\rightarrow \text{Mod-}\mathcal{H}(\tilde{G},\tilde{\lambda}),\,\qquad\tilde{\tau}\mapsto \hat{\tilde{\tau}}:=\Hom_{\tilde{J}}(\tilde{\lambda},\tilde{\tau}),$$
where the $\mathcal{H}_{\tilde{G}}$-action on $\hat{\tilde{\tau}}$ is given by
\begin{equation*}
(\phi\cdot f)(w)=\int_{\tilde{G}}\tilde{\tau}(g^{-1})\circ \phi\circ f(g)(w) dg,\qquad 
\phi\in \hat{\tilde{\tau}},\,f\in\mathcal{H}_{\tilde{G}},\,w\in V_{\tilde{\lambda}}.
\end{equation*}
(Here $V_{\tilde{\lambda}}$ is the representation space of ${\tilde{\lambda}}$, and recall that we fixed the Haar measure on $\tilde{G}$ such that the measure of $\tilde{J}_\Lambda$ is 1.) In particular, we have 
\begin{equation}
\label{shifted action of Z on eigenvector}
\phi\cdot Z= \tilde{\tau}(\varpi_E^{-1})\circ \phi\circ Z(\varpi_E).
\end{equation}

The group of unramified characters of $\tilde{G}$ acts on $\mathcal{R}^{[\tilde{G},\tilde{\pi}]}(\tilde{G})$ by 
 \begin{equation*}
 \tilde{\chi}\cdot \tilde{\tau} = \tilde{\tau}\otimes \tilde{\chi}.
 \end{equation*}
 and on $\mathcal{H}_{\tilde{G}}$ by $$(\tilde{\chi}\cdot f)(x)=\tilde{\chi}(x)f(x),\qquad f\in\mathcal{H}_{\Tilde{G}},\,x\in \tilde{G},$$ which induces an action on $\text{Mod-}\mathcal{H}_{\tilde{G}}$ naturally. The equivalence map $\mathcal{M}_{\tilde{G}}$ is hence equivariant under this action.

We now turn to our constructed type $\lambda_P$ in $G_W$ and abbreviate $\mathcal{H}_M=\mathcal{H}(M,\lambda_P|_{J_M})$ and $\mathcal{H}_{G_W}=\mathcal{H}(G_W,\lambda_P)$. For $s\in \mathbb{C}$, we denote by $D_s$ the simple right $\mathcal{H}_M$-module 
$$\mathcal{M}_M(\tilde{\pi}|\det|^s\boxtimes\pi)
=\Hom_{J_M}(\lambda_M,\tilde{\pi}|\det|^s\boxtimes\pi),$$ 
necessarily 1-dimensional. We then denote by $X_s$ be the right 
$\mathcal{H}_G$-module
$$\mathcal{M}(\tilde{\pi}|\det|^s\rtimes\pi)=\Hom_{J_P}(\lambda_P,\tilde{\pi}|\det|^s\rtimes\pi)\cong \Hom_{\mathcal{H}_M}(\mathcal{H}_G,D_s),$$
where the last isomorphism is given by the commutative diagram (\ref{commutative diagram}).

 Suppose that we choose $t_P$ normalized by the relation 
\begin{equation}
\label{embedding characterized by Z}
t_P(Z)(\zeta)=\Delta_P(\zeta)^{-1/2}Z(\varpi_E),
\end{equation}
where $\Delta_P$ is the modular character of $P$ and  
$\zeta =s_ys_z=i_M(\varpi_EI,I)$ as in the previous subsection, 
then by combining (\ref{shifted action of Z on eigenvector}), (\ref{embedding characterized by Z}), and (\ref{product Ty and Tz}) (evaluating at $\varpi_E$), the action of $t_P(Z)$ on $X_s$ is scalar, given by 
\begin{equation}
\label{product of Hecke algebra and parabolic induction}
q_E^{s}\Delta_P(\zeta)^{-1/2}{\tilde{\boldsymbol{\lambda}}}(\varpi_{E})
   T_y(s_y)T_z(s_z).
\end{equation}

We now impose an assumption that 
\begin{itemize}
    \item the operators $\tilde{\boldsymbol{\lambda}}(\varpi_{E})$, 
   $T_y(s_y)$, and $T_z(s_z)$ on $V_{\tilde{\lambda}}$ have finite orders.
\end{itemize}
This assumption can be easily satisfied when $\tilde{\boldsymbol{\lambda}}$ is self-dual and $s_y$ and $s_z$ are chosen to be simple enough; for example, we can and do pick \begin{equation}
\label{sy and sz}
s_y=\begin{bmatrix}
&&I_{V}
\\
&I_V&
\\
I_{V}&&
\end{bmatrix}
\quad \text{ and }\quad
s_z=\begin{bmatrix}
  &&-\varpi_E^{-1} I_{V}
  \\
  &I_{V}
&
  \\
\varpi_E I_{V}&&
\end{bmatrix},
\end{equation}
which by (\ref{T(s^2)}) means that $T_y(s_y)^2=1$ and $T_z(s_z)^2=\tilde{\rho}(-1)$. We also
use the fact 
\begin{equation*}  \Delta_P(\zeta)=[J_P^+: \zeta J_P^+\zeta^{-1}]=[s_y J_P^+s_y:J_P^-]^{-1}[s_z J_P^-s_z^{-1}:J_P^+]^{-1}, \end{equation*}
 which is just  $  (c_yc_z)^{-1}$.

We now compare the two quadratic relations (\ref{quadratic relation general b and c}) and (\ref{normalized quadratic relation general b and c}), and obtain the value of $b_w$, for $w\in \{y,z\}$, 
\begin{equation*}
 b_w = \pm c_w^{1/2} (q_E^{r_w/2  } - q_E^{ - r_w/2  }) .  
\end{equation*}
Hence the eigenvalues of $T_w$ are 
\begin{equation*}
    \pm \{c_w^{1/2} q_E^{r_w/2  }, -c_w^{1/2} q_E^{-r_w/2  }\}
\end{equation*}
and the possible products of eigenvalues of $T_y$ and those of $T_z$ are  
\begin{equation}
\label{product of eigenvalues of Ty and Tz}
    \pm (c_yc_z)^{1/2}\cdot \{ q_E^{\pm (r_y+r_z)/2  }, - q_E^{\pm (r_y-r_z)/2  }\}.
\end{equation}
When $\tilde{\pi}|\det|^s\rtimes\pi$, and hence $X_s$, is reducible, the eigenvalue of $Z$ is a product of one 
of those of $T_y$  
and one of those of $T_z$. By comparing the absolute values of (\ref{product of Hecke algebra and parabolic induction}) and (\ref{product of eigenvalues of Ty and Tz}), we obtain
\begin{equation}
\label{values of real parts}
    \{s_1,s_2\} = \{\frac{r_y+r_z}{2},\,\frac{|r_y-r_z|}{2}\}.
\end{equation}
In Section \ref{section First coefficient}, we will provide more detail about the products of the eigenvalues, and obtain the conditions on the cuspidal types for which reducibilities happen.

\subsection{M{\oe}glin's results}

In this subsection, we provide two results due to M{\oe}glin, for which we require the characteristic of $F$ to be 0. The first one is called a finiteness result in \cite[4.Prop]{Moeg-base-change}, and is improved in  \cite[Th.3.2.1]{Moeg-endosc-L-param} . The second one is a parity result on Langlands parameters \cite[5.6.Prop]{Moeg-base-change}. See also \cite[8.3.5]{Moeglin-Renard-non-quasi-split} for non-quasi-split groups.

 Given a supercuspidal representation $\pi$ of a unitary group $G=G_V$, then 
\begin{equation}
\label{Moeglin-estimate}
    \sum_{\tilde{\pi}}\sum_{ 
    \begin{smallmatrix}
    b\geq 0
    \\
b\equiv {2s_{\tilde{\pi}} -1} \mod 2
\end{smallmatrix}}
  b\dim V_{\tilde{\pi}} =  \dim V\end{equation}
where the sum $\sum_{\tilde{\pi}}$ ranges over all supercuspidal representations $\tilde{\pi}$, each of which is a supercuspidal representation of a general linear group $\tilde{G}_{V_{\tilde{\pi}}}$ for some space $V_{\tilde{\pi}}$,  such that $\tilde{\pi}|\det|^s\rtimes {\pi} $ is reducible at $s_{\tilde{\pi}}\in \frac{1}{2}\mathbb{Z}$ with $s_{\tilde{\pi}}\geq 1$. In particular, in the extreme case (which is studied in this paper), if there exists such a $\tilde{\pi}$ with $\dim V_{\tilde{\pi}} = \dim V$, then it is the unique such representation that gives rise to the reducibility above with $s_{\tilde{\pi}}=1$, in which case $\tilde{\pi}$ is the base change of ${\pi}_{} $.

We now switch to the Galois side and look at Langlands parameters. 
We recall that  $c$ is the generator of $\Gal({F/\Fo})$.       
We call an (semi-simple, smooth) irreducible representation $\Tilde{\varphi}$ of $\mathcal{W}_F$ conjugate-self-dual if $\tilde{\varphi}\cong {}^\sigma\tilde{\varphi}:= {}^c\tilde{\varphi}^\vee$ where
$${}^c\tilde{\varphi}(w)=\tilde{\varphi}(c^{-1}wc),\qquad 
\text{ for all }w\in \mathcal{W}_F,$$
and $\tilde{\varphi}^\vee$ is the contragredient of $\tilde{\varphi}$. Hence an irreducible representation
 $\tilde{\varphi}$ is conjugate-self-dual if there exists a non-degenerate bilinear form $B$ on $\mathbf{V}:=\mathbb{C}^{\deg \tilde{\varphi}}$ such that,
\begin{equation*} \label{B-is-G-invar}
  B({}^c\tilde{\varphi}(w)u,\tilde{\varphi}(w)v)=B(u,v)
  \qquad\text{ for all }u,v\in \mathbf{V}\text{ and }w\in \mathcal{W}_F.
\end{equation*}
We also define a parity on $\tilde{\varphi}$: call $\tilde{\varphi}$ conjugate-orthogonal (resp. conjugate-symplectic) if, furthermore,
\begin{equation}\label{B-is-Herm}
  B(u,v)= \sgn(\tilde{\varphi}) B(\tilde{\varphi}(c^2)v,u)
  \qquad\text{for all $u,v\in \mathbf{V}$},
\end{equation}
where $\sgn(\tilde{\varphi})=1$ (resp. $-1$). The parity of $\tilde{\varphi}$, as an irreducible representation of $\mathcal{W}_F$, can also be described using Asai L-functions: if $r_A$ is the Asai representation of $\GL(\mathbf{V})$ \cite[A.2.1]{Moeg-exhaustion}, then $\tilde{\varphi}$ is conjugate-orthogonal if and only if 
\begin{equation}
\label{Asai-L-for-parameter}
    L(s,\tilde{\varphi},r_A):= L(s,r_A\circ \tilde{\varphi}),\qquad s\in \mathbb{C},
\end{equation}
has a pole at $s=0$.

\section{Ramified unitary groups}
\label{section Ramified unitary group}

We now specify $F/\Fo$ to be a quadratic ramified extension, generated by a uniformizer $\varpi_F$ such that $\varpi_F^2=-N_{F/\Fo}(\varpi_F)=\varpi_\Fo.$

\subsection{The strongly ramified case}
\label{section Strongly ramified case}

Given a simple or null stratum $\mathbf{s}=[\Lambda,r,0,\beta]$, which is assumed to be skew, we assume the following condition on the field extension $E/\Eo$, where ${E}_\bullet$ is the subfield of $E$ fixed by $c = -\alpha$.   
 \begin{equation}
 \label{strong ramified condition ramified}
     \text{$E/E_\bullet $ is a quadratic ramified extension.}
 \end{equation}
The reason of imposing this condition has been explained at the end of the introduction section. This condition implies that $e=e(E/F)$ is an odd integer, and if $K$ is an intermediate subfield between $E$ and $F$, and $K_\bullet$ is the fixed-field of $K$ by the involution $-\alpha$, then  $K/K_\bullet $ is quadratic ramified.

Together with the simple stratum $ \mathbf{s}$ is a sequence of approximating simple strata $\mathbf{s}_j:=[\tilde{\mathfrak{A}},r,r_{j},\gamma_{j+1}]$, for $j=0,\dots,d$, satisfying the conditions in \cite[(2.4.2)]{BK}; in particular,

\begin{enumerate}[(i)]

\item the simple stratum $\mathbf{s}_j$ is equivalent to $[\tilde{\mathfrak{A}},r,r_{j},\gamma_{j}]$;

\item if we denote by

\begin{itemize}
\item

$\tilde{{B}}_j $  the centralizer of $\gamma_j$ in $\tilde{{A}}$,

\item $\tilde{\mathfrak{B}}_j  = \tilde{B}_j \cap  \tilde{\mathfrak{A}}$, and

\item $s_j$  the tame corestriction of $A$ relative to $F[\gamma_j]/F$, 
\end{itemize}

then the derived stratum $[\tilde{\mathfrak{B}}_{j+1},r_{j},r_{j}-1,s_{j+1}(\gamma_j-\gamma_{j+1})]$ is equivalent to a simple stratum, say $[\tilde{\mathfrak{B}}_{j+1},r_{j},r_{j}-1,\delta_{j}]$. 
\end{enumerate}

The numbers $\{r_j\}_{j=0}^d$ are the critical exponents of $\mathbf{s}$. For the moment we do not require the definition, but the simplicity of $[\tilde{\mathfrak{B}}_{j+1},r_{j},r_{j}-1,\delta_{j}]$ implies that 
$v_{\tilde{\mathfrak{B}}_{j+1}}(\delta_j)=-r_j$. Note that if $\mathbf{s}$ is skew, we can apply \cite[(1.10)]{stevens-intertwining-supercusp} and do choose all $\mathbf{s}_j$ to be also skew. If we also choose the tame corestriction $s_j$ to be $\alpha$-equivariant, then we can also assume  
$[\tilde{\mathfrak{B}}_{j+1},r_j,r_j-1,\delta_j]$ to be skew.

\begin{prop}\label{all jumps are odd}
  When $E/\Eo$ is ramified, all $\{r_j\}_{j=1}^d$ are odd.
\end{prop}
\proof
Let $[\tilde{\mathfrak{B}}_{j+1},r_j,r_j-1,\delta_j]$ be the simple stratum equivalent to the derived stratum $[\tilde{\mathfrak{B}}_{j+1},r_j,r_j-1,s_{j+1}(\gamma_j-\gamma_{j+1})]$, 
which is assumed to be skew; in  particular, we have $\delta_j={}^\alpha \delta_j$. Hence $-\alpha$ defines an involution on the field $F[\delta_j]$   
whose restriction to $F$ is $c$, the Galois conjugation of $F/\Fo$. We denote its fixed field by $F[\delta_j]_\bullet$. By \cite[(2.2.8)]{BK} and since $e(E/F)$ is odd, $e(F[\delta_j]/F])$ is also odd. This implies that $F[\delta_j]/F[\delta_j]_\bullet$ is ramified. The $\alpha$-invariance of $\delta_j$ implies that $v_{F[\delta_j]}(\delta_j)$ is odd, and so is $r_j = -v_{F[\delta_j]}(\delta_j) e(E/F)/e(F[\delta_j]/F)$. 
\qed

This proposition facilitates our calculations remarkably. First of all, we immediately see that $\tilde{J}^1=\tilde{H}^1$ from their constructions \cite[Sec. 3.1]{BK}. We also have $\tilde{J}^1_{\mathfrak{M}^y}= \tilde{H}^1_{\mathfrak{M}^y}$ since $\mathfrak{M}^y$ is a dilation and shift of a lattice sequence equivalent to $\Lambda$, and as for $\tilde{J}^1_{\mathfrak{m}}= \tilde{H}^1_{\mathfrak{m}}$ we just note that the critical exponents of $(\beta,\beta,\beta)$ are $\{3r_j\}_{j=1}^d$. The Heisenberg representation containing a simple character $\tilde{\theta}$ is just $\tilde{\theta}$ itself, and so any beta-extension of $\tilde{\theta}$ is a character. In particular, the unique p-primary beta extension $\tilde{\kappa}_0$ is trivial on $\boldsymbol{\mu}_E$. Finally, we choose our covering type $\lambda_P = \kappa_P \otimes \rho_M$ such that $$\kappa_P= \kappa^y_P,$$
 since in this case $\kappa_P|_{J_M} = \tilde{\kappa}_0 \boxtimes \kappa_0$, i.e., it is p-primary.

If moreover $\tilde{\theta}$ is self-dual, which is assumed from now on, then the p-primary beta extension $\tilde{\kappa}_0$ is a self-dual character. We fix a uniformizer $\varpi_\Eo$ of $\Eo$, and choose another one $\varpi_E$ for $E$ satisfying 
\begin{equation*}
\label{varpi-E and varpi-Eo}
\varpi_E^2=-N_{E/\Eo}(\varpi_E)=-\varpi_\Eo.
\end{equation*} 

\begin{prop}
\label{unique bold kappa}
Suppose that $\tilde{\theta}$ is a self-dual simple character, and let     
$\tilde{\kappa}_0$ be its p-primary beta-extension. There is a unique representation (a character) $\tilde{\boldsymbol{\kappa}}_0$ of $\tilde{\mathbf{J}}$ characterized by the following conditions:
\begin{enumerate}[(i)]
\item $\tilde{\boldsymbol{\kappa}}_0|_{\tilde{{J}}}=\tilde{\kappa}_0$; \label{unique bold kappa extension}

\item $\tilde{\boldsymbol{\kappa}}_0(\varpi_F)=1$;
\label{unique bold kappa varpi-F}

\item $\tilde{\boldsymbol{\kappa}}_0(\varpi_E)^2=\tilde{\boldsymbol{\kappa}}_0(-1)$. 
\label{unique bold kappa varpi-E}
\end{enumerate}
Moreover, such an extension $\tilde{\boldsymbol{\kappa}}_0$ is self-dual, i.e., ${}^\sigma \tilde{\boldsymbol{\kappa}}_0=\tilde{\boldsymbol{\kappa}}_0$. 
\end{prop}
\proof
Clearly there is a unique representation of $F^\times \tilde{{J}}$ satisfying (\ref{unique bold kappa extension}) and (\ref{unique bold kappa varpi-F}). Since $\varpi_E^e\in F^\times U_E$, the value of $\tilde{\boldsymbol{\kappa}}_0(\varpi_E)^e$ is known, and so is $\tilde{\boldsymbol{\kappa}}_0(\varpi_E)$ because of (\ref{unique bold kappa varpi-E}) and that $e$ is odd. The last statement is also clear since ${}^\sigma\varpi_E=-\varpi_E^{-1}$.
\qed

We also take a self-dual level-zero type $\tilde{\boldsymbol{\rho}}$ of $\tilde{\boldsymbol{J}} = E^\times\tilde{J} $, i.e., it is an extension of a representation $\tilde{{\rho}}$ of $\tilde{J}$ inflated from an irreducible cuspidal representation $\tilde{J}/\tilde{J}^1 \cong \GL_{f_0}(\mathbf{k}_E)$, where $f_0=n/[E:F]$. Denote $E^{f_0}$ the unramified extension of $E$ of degree $f_0$ in $\tilde A$.      
 Associated to $\tilde{\boldsymbol{\rho}}$ is a level-zero character $\tilde{\xi}_{\tilde{\boldsymbol{\rho}}}$ of $E^{f_0
\times}$ which is $\Gal(E^{f_0}/E)$-regular 
(which is equivalent to say that $\tilde{\xi}_{\tilde{\boldsymbol{\rho}}}|_{U_{E^{f_0}}^\times}$  is $\Gal(E^{f_0}/E)$-regular), such that $\tilde{\boldsymbol{\rho}}$ is self-dual if and only if $\tilde{\xi}_{\tilde{\boldsymbol{\rho}}}$ is, in which case it means that either $f_0=1$ or else, by \cite[Lemmas 2.1 and 5.1]{MR-unitary}, the involution ${-\alpha}$ on $E$ extends to $E^{f_0}$ which is unramified over the the fixed field $E^{f_0}_{\bullet} = (E^{f_0})^{-\alpha}$, and in particular $f_0$ is even.

We postpone the discussion of the case $f_0>1$  to   the appendix \ref{section appendix} of this section and proceed to our main results under the condition $f_0=1$.  The reason of imposing this condition 
is accounted for in  Proposition \ref{same parity strongly ramified}:  if $f_0>1$,   the  two candidates for the base  change  do not have the same parity, therefore the decision between them  can be made by computing Asai L-functions,   as was briefly explained in the last paragraph of the introduction.

So we  
 look at the case when 
\begin{equation}
    \label{strong ramified condition f=1}
    f_0=\dim_F V /[E:F]=1,
\end{equation}
so that both $\tilde{\boldsymbol{\kappa}}_0$ and $\tilde{\boldsymbol{\rho}} = \tilde{\xi}_{\tilde{\boldsymbol{\rho}}}$ are characters. In this case, 
 $\tilde{\kappa}_0|_{\boldsymbol{\mu}_E}
  $
 is trivial and $\tilde{\boldsymbol{\kappa}}_0(\varpi_E)^2=1.
$
A self-dual extended maximal simple type $\tilde{\boldsymbol{\lambda}}$ 
is a twist of $\tilde{\boldsymbol{\kappa}}_0$ by a self-dual character $\tilde{\boldsymbol{\rho}}$ of $E^\times$ trivial on $U_E^1$, which implies that  $\tilde{\rho}$, defined on  ${\boldsymbol{\mu}_E}$, is at most quadratic and $\tilde{\boldsymbol{\rho}}(\varpi_E)^2 =  \tilde{\rho}(-1) $, 
so that $\tilde{\boldsymbol{\rho}}(\varpi_E)$ is a 4-th root of unity. We summarize the above conditions in the following terminology.
\begin{dfn}
\label{The strongly ramified condition}
We call the extended type $\tilde{\boldsymbol{\lambda}}$ for $\tilde{G}$, as well as its induced  
 supercuspidal representation 
$\tilde{\pi}$, {\em strongly ramified}
when both (\ref{strong ramified condition ramified}) and (\ref{strong ramified condition f=1}) are satisfied, and  similarly for $\lambda$ and $\pi$ for $G$.
\qed\end{dfn}

We now recall from \cite[Sec. 5]{lusz-classical} (or see \cite[Sec. 7.6]{Lust-Stevens}) the values of the integers $r_w$, with $w\in \{y,z\}$, that appear in the quadratic relations in (\ref{normalized quadratic relation general b and c}) at the end of subsection \ref{section Structures of Hecke algebras}. From our construction we have $$\mathsf{G}_y \cong \mathsf{O}_{3,\mathbf{k}_E}, 
\qquad \text{so that} \qquad
r_y = 1,
$$
and also 
$$\mathsf{G}_z \cong \mathsf{Sp}_{2,\mathbf{k}_E}\times \mathsf{O}_{1,\mathbf{k}_E}
, 
\qquad \text{so that} \qquad
r_z = 
\begin{cases}
1 & \text{ if $\tilde{\rho}$ is trivial,}
\\
0 & \text{ otherwise.}
\end{cases}
 $$
 From (\ref{values of real parts}) the positive real parts of the points of reducibility are 
 \begin{equation*}
 \{0,1\}\text{ if $\tilde{\rho}$ is trivial,}
 \qquad\text{ and }\qquad
 \text{$\{1/2,1/2\}$ otherwise.}
 \end{equation*}
There is a support-preserving algebra morphism from \cite[(7.3)]{stevens-supercusp}, 
$$\mathcal{H}(U_{\mathfrak{M}^w,E},\tilde{\rho}\boxtimes \rho)\hookrightarrow \mathcal{H}(G_W,\lambda_P)\qquad \text{for each }w\in \{y,z\}. $$
Hence if we view each generator $T_w$ as an endomorphism of the module $X_s$ corresponding
to the parabolic induction $\tilde{\pi}|\det|^s\rtimes \pi$, then the above results, combining with the quadratic relation (\ref{quadratic relation general b and c}), imply that 
\begin{equation}
\label{the number}
b_w = \epsilon_w(q_E-1)(c_w / q_E)^{1/2},
\end{equation}
where $\epsilon_w $ is a sign. We will provide the precise value for the sign $\epsilon_{w}$ in the next subsection, with calculation given in Section \ref{section First coefficient}.

\subsection{Reducibility results}
\label{section Reducibility results}

To state the main result on the eigenvalues for the quadratic relations, we impose the following assumptions which are used in Section \ref{section First coefficient}.
\begin{itemize}
    \item We only consider a strongly ramified supercuspidal representation $\pi$ of $G$, i.e., if $(\theta,\rho)$ is a pair consisting of a simple character and a level zero cuspidal representation defining $\pi$, then that $E/\Eo$ is ramified and $\dim V/[E:F]  = 1$.

    \item The extension $E/F$ is tamely ramified, which allows us to assume that $E_j = F[\gamma_{j}]$ is contained in $E_{j-1}=F[\gamma_{j-1}]$ for all $j$, forming a tower of intermediate extensions between $E$ and $F$.
\end{itemize}

We also take a supercuspidal representation $\tilde{\pi}$ of $\tilde{G}$ constructed by $(\tilde{\theta},\tilde{\boldsymbol{\rho}})$ satisfying the same conditions as $\pi$ above, with  ${\theta}=(\tilde{\theta}|_{H^1})^{1/2}$. The following theorem will be proven in Section \ref{section First coefficient}:

\begin{thm}
\label{eigenvalues of Ty and Tz}
Suppose that $\tilde{\pi}$ and $\pi$ satisfy the above conditions.\begin{enumerate}[(i)]
\item  The coefficient $b_y$ of the quadratic relation of $T_y$ is
\begin{equation*}
b_y=\tilde{\rho}(-2)\rho(-1)T_y(s_y)(q_E-1)(c_y / q_E)^{1/2}. 
\end{equation*}

\label{eigenvalues of Ty}

\item  The coefficient $b_z$   
 of the quadratic relation of $T_z$ is given as follows.

\begin{enumerate}[(a)]

\item When $\tilde{\rho}|_{\boldsymbol{\mu}_E}=\left(\frac{\cdot}{\boldsymbol\mu_E}\right)^{f(E/F)}$, then
$
b_z=0$.

\item When $\tilde{\rho}|_{\boldsymbol{\mu}_E}=\left(\frac{\cdot}{\boldsymbol\mu_E}\right)^{f(E/F)-1}$, then 
\begin{equation*}
b_z=\tilde{\rho}(-2)\epsilon_z^P(\varpi_E,\mathbf{s})T_z(s_z)(q_E-1)(c_z / q_E)^{1/2},
\end{equation*}
where $\epsilon_z^P(\varpi_E,\mathbf{s})$ is a sign, associated to a quadratic Gauss sum and determined by the choice of $\varpi_E$ and the simple stratum $\mathbf{s}$.
\end{enumerate}
\label{eigenvalues of Tz}
\end{enumerate}
\qed\end{thm}

For the purpose of base change, the case when $\tilde{\rho}|_{\boldsymbol{\mu}_E}=\left(\frac{\cdot}{\boldsymbol\mu_E}\right)^{f(E/F)}$ is unimportant.

\begin{cor}
\label{The two eigenvalues}
\begin{enumerate}[(i)]
    \item The two eigenvalues of $T_y$ are 
$$-\tilde{\rho}(-2)\rho(-1)T_y(s_y)(c_y / q_E)^{1/2}
\qquad\text{and}\qquad
 \tilde{\rho}(-2)\rho(-1)T_y(s_y)(c_y / q_E)^{1/2}q_E. 
$$ 

\item When $\tilde{\rho}|_{\boldsymbol{\mu}_E}=\left(\frac{\cdot}{\boldsymbol\mu_E}\right)^{f(E/F)-1}$, the two eigenvalues of $T_z$ are 
$$-\tilde{\rho}(-2)\epsilon_z^P(\varpi_E,\mathbf{s})T_z(s_z)(c_z / q_E)^{1/2}
\qquad\text{and}\qquad
 \tilde{\rho}(-2)\epsilon_z^P(\varpi_E,\mathbf{s})T_z(s_z)(c_z / q_E)^{1/2} q_E .
$$   
\end{enumerate}
\qed\end{cor}

We now choose $\tilde{\pi} = \tilde{\pi}(\tilde{\theta},\tilde{\boldsymbol{\rho}})$ with 
$$\tilde{{\rho}}|_{\boldsymbol{\mu}_E}=\left(\frac{\cdot}{\boldsymbol\mu_E}\right)^{f(E/F)-1}
\quad\text{and}\quad
\tilde{\boldsymbol{\rho}}(\varpi_E)=\epsilon_z^P(\varpi_E,\mathbf{s}){\rho}(-1).$$

\begin{cor}
\label{Corollary points of reducibility}
The points of reducibility of $\tilde{\pi}|\det|^s\rtimes \pi$ are
\begin{equation*}
\text{$ \pm 1
\quad $ and $ \quad 
\frac{\pi i}{\log q_E}$.}
\end{equation*}
\end{cor}
\proof
Since the eigenvalue of $t_P(Z)$ is the product of eigenvalues of $T_y$ and $T_z$, the comparison in (\ref{product of Hecke algebra and parabolic induction}) gives
\begin{equation}
    \label{relation between lambda and the signs epsilon-y and epsilon-y}
    \tilde{\boldsymbol{\rho}}(\varpi_E)q_E^s
=\epsilon_z^P(\varpi_E,\mathbf{s}){\rho}(-1)\cdot(-1 \text{ or } q_E^{\pm 1} ).
\end{equation}
The corollary follows by solving $s$.\qed 

\begin{rmk}
If $\tilde{\rho}|_{\boldsymbol{\mu}_E}=\left(\frac{\cdot}{\boldsymbol\mu_E}\right)^{f(E/F)}$, then we 
put $\tilde{\boldsymbol{\rho}}(\varpi_E)=\pm{\rho}(-1)\tilde{\rho}(-2)\tilde{\rho}(-1)^{1/2}$ and use similar arguments as above to show that points of reducibility are
\begin{equation*}
\text{$\pm \frac{1}{2} 
\quad $ and $ \quad 
\pm \frac{1}{2}+\frac{\pi i}{\log q_E}$}.
\end{equation*}
The choice of the square root $\tilde{\rho}(-1)^{1/2}$ for $\tilde{\boldsymbol{\rho}}(\varpi_E)$ does not matter as we will discard this case for base change anyway. \qed\end{rmk}

The following remark explains Corollaries \ref{The two eigenvalues} and \ref{Corollary points of reducibility} are independent of the various choices made throughout the progress. 

\begin{rmk}
Note that both $T_w(s_w)$, for $w\in \{y,z\}$, have been cancelled out in (\ref{relation between lambda and the signs epsilon-y and epsilon-y}). Indeed the reducibility points, as well as the base change result in Theorem \ref{main theorem} below, are independent of the normalizations of $T_w$ (and also that of $Z$) and the choices of $s_w$. Moreover, as explained in Subsection \ref{section Structures of Hecke algebras}, if we have chosen $P^-$ instead of $P$, then we need to switch the roles of $y$ and $z$. In this case, $\epsilon_z^P(\varpi_E,\mathbf{s})$ should then be denoted by $\epsilon_y^{P^-}(\varpi_E,\mathbf{s})$. Our corollaries are clearly independent of choosing $P$ or $P^-$. 
\qed\end{rmk}

\subsection{The main result for base change}
\label{section main result}

Let's first look at a simple case when the representations are characters, in which case $E=F$. The base change $\Tilde{\rho}$ of a character $\rho$ from $\mathrm{U}_1(F/\Fo)$ to $\GL_1(F)$ is $$\Tilde{\rho}(x) = \rho(x{}^\sigma x),
\qquad\text{for }x\in F^\times$$
or more explicitly, 
\begin{equation}
\label{explicit base-change for U1}
 {\rho}|_{(U_F^1)^\sigma}=(\tilde{\rho}|_{(U_F^1)^\sigma})^{1/2},\qquad \tilde{\rho}|_{\boldsymbol{\mu}_F}\equiv 1,
 \qquad\text{ and }
 \qquad
 \tilde{\rho}({\varpi_F})=\rho(-1).
\end{equation}
 Note that the restriction $\Tilde{\rho}|_{\boldsymbol{\mu}_F}$ has no relation with the level-zero part of $\rho$, i.e., the character of the finite reductive quotient of $\mathrm{U}_1$, which is $\mathrm{O}_1= \{\pm 1\}$. However, if we change $\Tilde{\rho}$ by twisting a quadratic unramified character, then  $\rho(-1)$ is changed to another sign. 

We now consider strongly ramified supercuspidal representations  $\pi=\pi(\theta,\rho)$ and $\tilde{\pi} = \tilde{\pi}(\Tilde{\theta},\Tilde{\boldsymbol{\rho}})$ as in Definition \ref{The strongly ramified condition}. The following theorem, the main result of our paper, provides the relations between $(\theta,\rho)$ and $(\Tilde{\theta},\Tilde{\boldsymbol{\rho}})$ for which $\tilde{\pi}$ is the base change of $\pi$.%
\begin{thm}
\label{main theorem}
Suppose that char$(F)=0$, and the simple characters ${\theta}$ and $\tilde{\theta}$ and the tamely ramified characters $\tilde{\boldsymbol{\rho}}$ and $\rho$ are related as follows.
\begin{enumerate}[(i)]
\item ${\theta}=(\tilde{\theta}|_{H^1})^{1/2}$,

\item $\tilde{\rho}|_{\boldsymbol{\mu}_E}=\left(\frac{\cdot}{\boldsymbol{\mu}_E}\right)^{f(E/F)-1}$, and

\item 
$\tilde{\boldsymbol{\rho}}(\varpi_E)=\rho(-1)\epsilon_z^P(\varpi_E,\mathbf{s})$, where $\epsilon_z^P(\varpi_E,\mathbf{s})$ is the sign appearing in Theorem \ref{eigenvalues of Ty and Tz}.
\end{enumerate}
Then $\tilde{\pi}$ is the base change of $\pi$. 
\qed\end{thm}

This can be easily deduced from the reducibility result in the previous subsection, together with the finiteness result of Moeglin (\ref{Moeglin-estimate}). One may notice that  (\ref{explicit base-change for U1}) is a special case of the theorem.

\begin{rmk}
\cite[7.1]{Moeg-base-change} implies that if in general a discrete series parameter, when viewed as a representation of the Weil-Deligne group, has $k$ irreducible components, then its corresponding L-packet contains $2^{k-1}$ members. In our case when the base change representation is supercuspidal, then its parameter is an irreducible representation, and the L-packet is a singleton. 
\qed\end{rmk}

\subsection{Appendix: the strongly ramified case and parity}
\label{section appendix}

This appendix is a sequel to subsection    
\ref{section Strongly ramified case}, and does not intervene with our main results.

We provide a result on the parity of a conjugate self-dual supercuspidal representation. To this end, we have to switch to the Galois side via the local Langlands correspondence for $\GL_n$, and assume that char$(F)=0$. We call a conjugate self-dual supercuspidal representation $\tilde{\pi}$ of a general linear group conjugate-orthogonal (resp. conjugate-symplectic) if the Asai L-function \cite{shahidi-complementary} for $\tilde{\pi}$, 
\begin{equation*}
    L(s,\tilde{\pi},r_A),\qquad s\in \mathbb{C},
\end{equation*}
has a pole at $s=0$. By \cite{Hen-ext-sym}, we know that $\tilde{\pi}$ is conjugate-orthogonal (resp. conjugate-symplectic) if its Langlands parameter is so (see {\ref{B-is-Herm}) and (\ref{Asai-L-for-parameter})).

\begin{rmk}
A conjugate self-dual character $\tilde{\chi}$ of $\GL_1(F)$ (i.e., ${}^\sigma \tilde{\chi} = \tilde{\chi}$) is  conjugate-orthogonal (resp. conjugate-symplectic) if and only if $\tilde{\chi}|_{\Fo^\times }$ is trivial (resp. is equal to $\delta_{F/\Fo}$, the character on $\Fo^\times$ with kernel the norm group $N_{F/\Fo}(F^\times)$). In particular, the quadratic unramified character $|\det|^{\pi i/\log q}$ is conjugate-orthogonal (resp. conjugate-symplectic) when $F/\Fo$ is ramified (resp. unramified).
\qed\end{rmk}

\begin{prop}
\label{same parity strongly ramified}
 When $F/\Fo$ is ramified, a conjugate self-dual supercuspidal representation  $\tilde{\pi}$ and its twist $\tilde{\pi}' := \tilde{\pi}|\det|^{\pi i/f_0 \log q_E}$ are of the same parity if and only if $\tilde{\pi}$ is strongly ramified. 
\end{prop}
\proof
Let $T/F$ be the maximal tamely ramified subextension of $E/F$, and $T^{f_0}$ be the unramified extension of $T$ of degree ${f_0}$ with fixed field $T^{f_0}_\bullet = T^{f_0}\cap E^{f_0}_\bullet$. Note that $T^{f_0}/T^{f_0}_\bullet$ is ramified if and only if $\tilde{\pi}$ is strongly ramified. From \cite[Ch.1]{bh-eff}, suppose the Langlands parameter of $\tilde{\pi}$ takes the form 
\begin{equation*}
    \Ind_{T/F}(\Tilde{\alpha}\otimes \Ind_{T^{f_0}/T}\Tilde{\xi}) \cong \Ind_{T^{f_0}/F}(\Res_{T^{f_0}/T}\Tilde{\alpha}\otimes \Tilde{\xi}_{}) 
\end{equation*}
where $\Tilde{\alpha}$ is an irreducible representation of $\mathcal{W}_{T}$, and $\tilde{\xi_{}}$, as a character of $T^{{f_0}\times}$, is viewed as a character of $\mathcal{W}_{T^{f_0}}$. The Langlands parameter of $\tilde{\pi}'$ is $$ \Ind_{T^{f_0}/F}(\Res_{T^{f_0}/T}\Tilde{\alpha}\otimes \Tilde{\xi}_{} |\mathrm{det}_{T^{f_0}}|^{\pi i/\log q_{T^{f_0}}}) .$$
Since the parity of a representation is preserved under induction, the result then holds by the remark before the proposition, and the fact that $|\mathrm{det}_{T^{f_0}}|^{\pi i/\log q_{T^{f_0}}}$ is conjugate-symplectic when ${f_0}$ is even.
\qed

\section{The coefficients }\label{section First coefficient}

The whole section is devoted to prove Theorem \ref{eigenvalues of Ty and Tz}. It suffices to compute the values of $b_w$ for $w\in\{y,z\}$. We do not require that char$(F)=0$.

We begin by recalling the explicit form of $J^1_P$. First of all, we can express the rings as follows 
\begin{equation*}
\label{expanding H and J}
\begin{split}
\tilde{\mathfrak{H}} & = \tilde{\mathfrak{A}}_{{\Lambda},E} + 
\tilde{\mathfrak{P}}^{(r_0/2)_+}_{{\Lambda},E_{1}}
+ \cdots 
+\tilde{\mathfrak{P}}^{(r_{d-1}/2)_+}_{{\Lambda},E_{d}} 
+ \tilde{\mathfrak{P}}^{(r_d/2)_+}_{{\Lambda}},
\\
\tilde{\mathfrak{J}} & = \tilde{\mathfrak{A}}_{{\Lambda},E} + 
\tilde{\mathfrak{P}}^{r_0/2}_{{\Lambda},E_{1}}
+ \cdots 
+\tilde{\mathfrak{P}}^{r_{d-1}/2}_{{\Lambda},E_{d}} 
+ \tilde{\mathfrak{P}}^{r_d/2}_{{\Lambda}}.
\end{split}
\end{equation*} 
 By \cite[Prop. 1]{blondel-propag},
\begin{equation}
\label{explicit form of J^1_m}
\mathfrak{J}_\mathfrak{m}^1=\begin{bmatrix}
\tilde{\mathfrak{J}}^1 & \tilde{\mathfrak{J}}^0& \varpi_E^{-1}\tilde{\mathfrak{H}}^1
\\
\tilde{\mathfrak{H}}^1 & \tilde{\mathfrak{J}}^1& \tilde{\mathfrak{J}}^0
\\
\varpi_E\tilde{\mathfrak{J}}^0 & \tilde{\mathfrak{H}}^1& \tilde{\mathfrak{J}}^1
\end{bmatrix}\cap A_W,
\end{equation}
In the strongly ramified case, since all exponents $r_i=2s_i+1$, for $i=0,\dots,d$, are odd, we have $\tilde{\mathfrak{J}}_\Lambda = \tilde{\mathfrak{H}}_\Lambda$. We can show that $\mathfrak{J}_\mathfrak{m}^1=\mathfrak{H}_\mathfrak{m}^1$ and ${J}_\mathfrak{m}^1={H}_\mathfrak{m}^1=J^1_P$.

Note that we have chosen certain convenient normalizations of $T_w$, for $w\in \{y,z\}$, and $Z$ in the relation $t_P(Z)=T_y*T_z$ to simplify our calculations. Following \cite{Bl-Bl-SP4}, we may choose $T_w$
 such that 
\begin{equation*}
T_w(s_w) T_w(s_w^{-1})=1,
\end{equation*}
which can be easily shown to be equivalent to requiring that $T_w(s_w)^2=\lambda_P(s_w^2)$.

\subsection{Computation of $b_y$}

We use (\ref{formula for b and c}) to compute $b_y$, so we compute  $T_y((X,Y)^-)$ for $(X,Y)^-\in {s_yJ_P^+s_y \cap J_Ps_yJ_P}/{J_P^-}$. We also recall (\ref{relation between X and Y}) which we use repeatedly: 
\begin{equation}
\label{relation between X and Y recalled}
X{}^\alpha X=Y-{}^\alpha Y.
\end{equation}
We take
\begin{equation*}
s_y=\begin{bmatrix}
&&I
\\
&I&
\\
I&&
\end{bmatrix}. 
\end{equation*}
If $(X,Y)^-=\left[\begin{smallmatrix}
I&&
\\
{}^\alpha X&I&
\\
Y&X&I
\end{smallmatrix}\right]\in {s_yJ_P^+s_y }/{J_P^-}$, then from (\ref{explicit form of J^1_m}) we have,
\begin{equation*}
X\in \tilde{\mathfrak{J}}^0/\tilde{\mathfrak{H}}^1
\quad\text{ and }\quad
Y\in\varpi_E^{-1}\tilde{\mathfrak{H}}^1/\varpi_E \tilde{\mathfrak{J}}^0. 
\end{equation*}
Also, we write 
$$\mathrm{supp} (T_y)= J_P  s_y J_P= J_P^+ J_M s_y J_P^+,$$ 
so that $(X,Y)^- \in s_yJ_P^+s_y \cap J_Ps_yJ_P$ can be written as
\begin{equation*}
\begin{split}
\begin{bmatrix}
   I&&  \\ {}^\alpha X &I& \\  Y&X&I
\end{bmatrix}
=&\begin{bmatrix}
   I&{}^\sigma YX&Y^{-1}  \\ &I&{}^\alpha({}^\sigma YX)\\  &&I
\end{bmatrix}
\begin{bmatrix}
   && I \\ &I& \\ I &&
\end{bmatrix}
\\
&\begin{bmatrix}
   Y&&  \\ &I-{}^\alpha XY^{-1} X& \\  &&{}^\sigma Y
\end{bmatrix}
\begin{bmatrix}
   1&Y^{-1}X &Y^{-1}  \\ &1&{}^\alpha (Y^{-1}X) \\  &&1
\end{bmatrix},
\end{split}
\end{equation*}
or in simplified symbols, 
\begin{equation}
\label{y-element-into-supp-product}
(X,Y)^-=({}^\sigma YX,Y^{-1})^+
\cdot  s_y 
\cdot i_M(Y,I-{}^\alpha XY^{-1} X)\cdot(Y^{-1}X,Y^{-1})^+.
\end{equation}

\begin{lem}
\label{representatives for X and Y for y}
Each coset in $(s_yJ_P^+s_y \cap  J_P s_y J_P) / J_P^-$ has a representative of the form 
$(X,Y_0(I+Y'))^-$ with $Y_0\in \boldsymbol{\mu}_E$ and $X\in( \mathfrak{o}_{E} \setminus \mathfrak{p}_{E} )\mod\mathfrak{p}_{E}$  such that $2Y_0\equiv -X^2 \mod {\mathfrak{p}}_E$, and $I+Y'\in \tilde{J}^1$. 
\end{lem}
\proof
The coset space $\tilde{\mathfrak{J}}^0/\tilde{\mathfrak{H}}^1$ containing $X$, when viewed as a $\mathbf{k}_F$-vector space, takes the form
$ (\mathfrak{o}_E/\mathfrak{p}_E)
\oplus (\tilde{\mathfrak{J}}^1/\tilde{\mathfrak{H}}^1)$.
If $r_j=2s_j+1$ are odd,
then $\tilde{\mathfrak{J}}^1= \tilde{\mathfrak{H}}^1$, and so we can choose $X\in \mathfrak{o}_{E} \mod\mathfrak{p}_{E}$. If furthermore $(X,Y)^-\in J_Ps_yJ_P$, then $Y $ is forced    to belong to $\tilde{J}$ by (\ref{y-element-into-supp-product}), which allows us to choose $Y_0 \in \boldsymbol{\mu}_E$. The relation (\ref{relation between X and Y recalled}) implies that $2Y_0\equiv -X^2 \mod \mathfrak{p}_E$, so that $X \notin \mathfrak{p}_E$, and $Y_0$ is uniquely determined by $X$. Hence the lemma follows. \qed

From (\ref{y-element-into-supp-product}) we obtain
\begin{equation*}
\label{t_y mid-calculation}
T_y((X,Y)^-)=  T_y(s_y) \tilde{\lambda}(Y)\lambda(I-{}^\alpha XY^{-1}X).
\end{equation*}
The following lemma shows that it is indeed a constant.

\begin{lem}
\label{T_y constant}
$
T_y((X,Y)^-)=\tilde{\rho}(-2)\rho(-1)T_y(s_y)
$
for all $(X,Y)^-\in (s_yJ_P^+s_y\cap J_Ps_yJ_P)/J_P^-$. 
\end{lem}
\proof
Since now $X\in  \tilde
{J}$ is invertible, we have 
\begin{equation*}
I-{}^\alpha XY^{-1}X = -{}^\alpha X Y^{-1}{}^\sigma Y^{-1}{}^\alpha X^{-1}
\end{equation*}
and so 
\begin{equation*}
\lambda(I-{}^\alpha XY^{-1}X) = 
\rho(-1)\theta(Y{}^\sigma Y)^{-1}
=
\rho(-1)\theta((I+Y'){}^\sigma(I+Y'))^{-1}.    
 \end{equation*}  
On the other hand, that $2Y_0\equiv -X^2 \mod {\mathfrak{p}}_E$ implies that 
\begin{equation*}
\tilde{\lambda}(Y) = \tilde{\rho}(Y_0)\tilde{\theta}(I+Y')=\tilde{\rho}(-2)\tilde{\theta}(I+Y').
\end{equation*}
Since $\tilde{\theta}$ and $\theta$ are related by (\ref{simple-char-base-change}), we have 
$$
\tilde{\theta}(I+Y')=\theta((I+Y'){}^\sigma(I+Y')),
$$
 and the lemma follows.
\qed

By combining the above two lemmas with (\ref{the number}), we obtain      
\begin{equation*}
b_y=\tilde{\rho}(-2)\rho(-1)T_y(s_y)(q_E-1)(c_y/q_E)^{1/2}.
\end{equation*}
This proves Theorem \ref{eigenvalues of Ty and Tz}.(i).

\subsection{Computation of $b_z$}
\label{section coefficient b_z}

We use (\ref{formula for b and c for y }) to compute $b_z$, so we compute  $T_z((X,Y)^+)$ for $(X,Y)^+\in {s_z^{-1}J_P^-s_z \cap J_P  s_z J_P}/{J_P^+}$. We take
$$s_z=\begin{bmatrix}
  &&-\varpi_E^{-1} I
  \\
  &I&
  \\
  \varpi_E I&&
\end{bmatrix}.$$
If $(X,Y)^+=\left[\begin{smallmatrix}
I&X&Y
\\
&I&{}^\alpha X
\\
&&I
\end{smallmatrix}\right]\in {s_z^{-1}J_P^-s_z }/{J_P^+}$, then we have
\begin{equation*}
X\in \varpi_E^{-1}\tilde{\mathfrak{H}}^1/\tilde{\mathfrak{J}}^0
\quad\text{ and }\quad
Y\in \varpi_E^{-1}\tilde{\mathfrak{J}}^0/\varpi_E^{-1}\tilde{\mathfrak{H}}^1. 
\end{equation*}
We also write
$$\mathrm{supp} (T_z)= J_P  s_z J_P= J_P^-   J_M s_z J_P^-,$$ 
so that $(X,Y)^+\in J_Ps_zJ_P \cap s_zJ_P^-s_z $ can be written as
\begin{equation*}
\begin{split}
\begin{bmatrix}
  I&X&Y  \\ &I&{}^\alpha X \\  &&I
\end{bmatrix}
=&\begin{bmatrix}
   I&& \\ {}^\alpha XY^{-1}&I& \\  Y^{-1} &{}^\alpha ({}^\alpha XY^{-1})&I
\end{bmatrix}
\begin{bmatrix}
   -Y\varpi_E&&  \\ &I-{}^\alpha XY^{-1} X& \\  &&{}^\sigma Y\varpi_E^{-1}
\end{bmatrix}
\\
&\begin{bmatrix}
   && -\varpi_E^{-1}I \\ &I& \\ \varpi_EI &&
\end{bmatrix}
\begin{bmatrix}
  I&&  \\ {}^\alpha X{}^\sigma Y  &I& \\  Y^{-1}&{}^\alpha ({}^\alpha X{}^\sigma Y )&I
\end{bmatrix},
\end{split}
\end{equation*}
or in simplified symbols, 
\begin{equation}\label{z-element-into-supp-product}
(X,Y)^+=({}^\sigma YX,Y^{-1})^-
\cdot i_M( -Y\varpi_E,I-{}^\alpha XY^{-1} X)
\cdot s_z
\cdot(Y^{-1}X,Y^{-1})^-.
\end{equation}
If we write $
Y= y\varpi_E^{-1}(I+Y')$ with $y\in \mathfrak{o}_E \mod \mathfrak{p}_E$ and $I+Y'\in \tilde{J}^1 = \tilde{H}^1$, then from (\ref{z-element-into-supp-product}) we can assume that $y\in \boldsymbol{\mu}_E$. We hence obtain
\begin{equation}\label{TzXYplus}
T_z((X,Y)^+)=\tilde{\rho}(-y)\tilde{\theta}(I+Y')\theta(I-{}^\alpha XY^{-1} X) T_z(s_z).
\end{equation}
In the following subsections, we will expand and simplify $\tilde{\theta}(I+Y')$ and $\theta(I-{}^\alpha XY^{-1} X)$.

\subsubsection{Expanding simple characters}

We fix an additive character $\psi$ of $F$ of conductor 1, which means that it is trivial on $\mathfrak{p}_F$ but not $\mathfrak{o}_F$. If $c\in A$ with $v_\Lambda(c)=-r$, then 
\begin{equation*}
\psi_{c}:1+X\mapsto \psi_F\circ\tr_{A/F}(c X)
\end{equation*}
defines a character on the compact subgroup $\tilde{U}^{(r/2)_+}_{\Lambda} $ which looks additive, in the sense that 
\begin{equation*}
\psi_{c}((I+X)(I+Y))=\psi_{c}(I+X+Y).
\end{equation*}

Recall from section \ref{section Strongly ramified case}      
that, given the skew simple stratum $\mathbf{s}=[\Lambda,r,0,\beta]$, we have an approximation by skew simple strata $[\Lambda,r,r_{j},\gamma_{j+1}]$, for $j=0,\dots,d$, equivalent to $[\Lambda,r,r_{j},\gamma_{j}]$. We denote $c_j=\gamma_{j}-\gamma_{j+1}$ (where we take $\gamma_{d+1}=0$) and $E_j=F[\gamma_j]$.
We have $v_\Lambda(c_j) = -r_j$. Since $r_j=2s_j+1$ is odd, $(r_j/2)_+=s_j+1$.

The compact subgroup $\tilde{H}^1$ factorizes as (by \cite[(3.1.15)]{BK} inductively) 
\begin{equation}
\label{H^1-tilde factorization}
\tilde{U}_{E}^1
\tilde{U}^{s_0+1}_{\Lambda,{E_1}}  \cdots \tilde{U}^{s_{d-1}+1}_{\Lambda,{E_d}} 
\tilde{U}^{s_d+1}_{\Lambda} .
\end{equation}
Denote $\tilde{H}^{t}=\tilde{H}^1\cap \tilde{U}^t_{\Lambda} $ for all $t\geq 1$. A simple character $\tilde{\theta}=\tilde{\theta}_0\in  \tilde{\mathcal{C}}(\Lambda,0,
\beta)$ takes the following inductive form: for $j=0,\dots,d+1$, 
\begin{equation*}
\tilde{\theta}_j|_{ \tilde{U}^{s_{j-1}+1}_{\Lambda,{E_j}} }
=\tilde{\xi}_j\circ\mathrm{det}_{\tilde{{B}}_j}
\end{equation*}
for some character $\tilde{\xi}_j$ of $U^{s_{j-1}+1}_{E_j}$, and for $j=0,\dots,d$,
\begin{equation*}
\tilde{\theta}_j|_{\tilde{H}^{s_{j}+1}_{} }=\tilde{\theta}_{j+1}\psi_{c_{j}} 
\end{equation*}
for some simple character $\tilde{\theta}_{j+1}\in \tilde{\mathcal{C}}(\Lambda,s_{j},
\gamma_{j})$.

Similarly, $\theta$ is defined on 
\begin{equation}
\label{H^1 factorization}
H^1={U}^1_{{E}/\Eo}
{U}^{s_0+1}_{\Lambda,{{E_1}/\Eo_1}} \cdots {U}^{s_{d-1}+1}_{\Lambda,{{E_d}/\Eo_d}}
{U}^{s_{d}+1}_{\Lambda,{{F}/\Fo}}
\end{equation}
and    
takes the form
\begin{equation*}
{\theta}_j|_{ {U}^{s_{j-1}+1}_{\Lambda,E_j/\Eo_j} }
={\xi}_j\circ\mathrm{det}_{\tilde{{B}}_j}\
\end{equation*}
for some character ${\xi}_j$ of $U^{s_{j-1}+1}_{E_j/\Eo_j}$ and $j=0,\dots,d+1$, and
\begin{equation*}
{\theta}_j|_{{H}^{s_{j}+1}_{} }={\theta}_{j+1}\psi_{c_{j}/2}
\end{equation*}
for some simple character $ {\theta}_{j+1}\in \mathcal{C}(\Lambda,s_{j},
\gamma_{j})$ and $j=0,\dots,d$.

Since $\tilde{\theta}$ and $\theta$ are related by (\ref{simple-char-base-change}), we can and do assume similar relations between $\tilde{\theta}_j$ and $\theta_j$, and also between $\tilde{\xi}_j$ and $\xi_j$, i.e., 
\begin{equation*}
\xi_j=(\tilde{\xi}_j|_{U^{s_{j-1}+1}_{E_j/\Eo_j}})^{1/2}.
\end{equation*}

In the following subsections, we call the factors of $\tilde{\theta}$ involving $\tilde{\xi}_j$ the multiplicative parts of  $\tilde{\theta}$, and those involving $\psi_{c_j}$ the additive parts, and similarly for $\theta$.

\subsubsection{Cancellation of multiplicative parts}

We recall that we have written $Y= y\varpi_E^{-1}(I+Y')$, with $y\in \boldsymbol{\mu}_E$ and $I+Y' \in \tilde{H}^1$.    
We first compare the multiplicative parts of $\tilde{\theta}(I+Y')$ and $\theta(I-{}^\alpha XY^{-1} X)$ by rewriting 
\begin{equation*}
    \tilde{\theta}(I+Y') = {\theta}((I+Y'){}^\sigma (I+Y')) = \rho(-1)\theta(I-X{}^\alpha X Y^{-1})^{-1},
\end{equation*}
and so we are actually comparing 
\begin{equation*}
    \theta(I-{}^\alpha XY^{-1} X)
    \qquad\text{and}\qquad 
    \theta(I-X{}^\alpha X Y^{-1}).
\end{equation*}

Starting from the proposition below, we have to assume that 
\begin{equation*}
\text{$E/F$ is tamely ramified.}
\end{equation*}
This condition allows us to assume that each field $E_j=F[\gamma_j]$, is contained in $E_{j-1}$ for all $j$, forming a tower of intermediate extensions between $E$ and $F$.

\begin{lem}\label{propmultcancel}
The multiplicative parts of $ \theta(I-{}^\alpha XY^{-1} X)
   $ and $\tilde{\theta}(I+Y')$ cancel with each other.
\end{lem}

\proof
Let's temporarily write $W={}^\alpha X Y^{-1}$. We have to compare $ \theta(I-W X)$ and $
    \theta(I-XW)$.
    Let's first expand $X = X_0+\cdots +X_{d+1}$, and similarly for $W$. For each $i=0,\dots,d+1$, we write $X_{\leq i} = X_0 + \cdots +X_i$, and similarly for $W_{\leq i}$. The part of the character $\theta(I-WX )$ involving ${\xi}_i$ is 
    \begin{equation*}
        {\xi}_i\circ\det{}_{E_i}((I-W_{\leq i-1}X_{\leq i-1})^{-1}
        (I-W_{\leq i}X_{\leq i}))
    \end{equation*}
    and similarly for the part of $  \theta(I-X{}^\alpha X Y^{-1})$ involving ${\xi}_i$, with $W$ and $X$ exchanged. Due to the identity 
\begin{equation*}
    \det{} (I-W X ) = \det{} (I-X W)
 \end{equation*}
 these two parts are the same. 
\qed

\subsubsection{Choosing representatives}

In this subsection, we will expand $X$ and $Y= y \varpi_E^{-1}(I+Y')$ such that the additive parts of $\tilde{\theta}(I+Y')$ and of $\theta(I-{}^\alpha XY^{-1}X)$ admit many simplifications.

The coset spaces $\varpi_E^{-1}\tilde{\mathfrak{H}}^1 / \tilde{\mathfrak{J}}^0$ and $\varpi_E^{-1}\tilde{\mathfrak{J}}^0 /\varpi_E^{-1}\tilde{\mathfrak{H}}^1 
$ respectively containing $X$ and $Y$, when viewed as $\mathbf{k}_F$-vector spaces, take the form
\begin{equation}
\label{coset space W-z}
\varpi_E^{-1}\tilde{\mathfrak{H}}^1 / \tilde{\mathfrak{J}}^0\cong  \bigoplus_{j=0}^d\mathfrak{W}_{z,j},\quad\text{where}\quad
\mathfrak{W}_{z,j}= \tilde{\mathfrak{P}}^{-1+(r_j/2)_+}_{\Lambda,E_{j+1}}/ (\tilde{\mathfrak{P}}^{-1+(r_j/2)_+}_{\Lambda,E_{j}}+ \tilde{\mathfrak{P}}^{r_j/2}_{\Lambda,E_{j+1}})
\end{equation}
and
\begin{equation*}
\varpi_E^{-1}\tilde{\mathfrak{J}}^0 /\varpi_E^{-1}\tilde{\mathfrak{H}}^1 \cong \mathbf{k}'_E\oplus  \bigoplus_{j=0}^d\mathfrak{W}'_{z,j},
\end{equation*}
where $\mathbf{k}'_E\cong \mathfrak{p}_E^{-1}/\mathfrak{o}_E$ and $
\mathfrak{W}'_{z,j}= \tilde{\mathfrak{P}}^{-1+(r_j/2)}_{\Lambda,E_{j+1}}/ (\tilde{\mathfrak{P}}^{-1+(r_j/2)}_{\Lambda,E_{j}}+ \tilde{\mathfrak{P}}^{-1+(r_j/2)_+}_{\Lambda,E_{j+1}})$. Since each $r_j=2s_j+1$ is odd, the summand $\mathfrak{W}'_{z,j}$ is trivial. We expand $X=\sum_{j=0}^{d+1}X_j$ and $Y'=\sum_{j=0}^{d+1}Y'_{j}$ accordingly, first requiring that
\begin{equation*}
Y_0 =  y\varpi_E^{-1} \text{ where }y \in \boldsymbol{\mu}_E,
\quad\text{and}\quad 
X_j  \in \tilde{\mathfrak{P}}^{s_{j-1}}_{\Lambda,{E_j}}  \mod 
\tilde{\mathfrak{P}}^{s_{j-1}+1}_{\Lambda,{E_j}}  + \tilde{\mathfrak{P}}^{s_{j-1}}_{\Lambda,{E_{j-1}}} \text{ for }j>0.
\end{equation*}
With these fixed, we choose auxiliary data 
\begin{equation}
\label{Groth-group for X in z}
X_0\in \mathfrak{o}_E
\quad\text{and}\quad 
Y_j \in \tilde{\mathfrak{P}}^{s_{j-1}}_{\Lambda,{E_j}}\text{ for }j>0
\end{equation}
 such that 
$
X{}^\alpha X = Y - {}^\alpha Y$ still holds. Eventually our main results are independent of these auxiliary choices, see (\ref{additivepartfinal}) for example.

We require some notations. Let $E_i = F[\gamma_i]$ and $ \tilde{B}_i$ the centralizer of $\gamma_i$ in $\tilde{A}$. For $i\geq 1$, we denote by 
\begin{equation*}
\text{$\tilde{B}_{i-1}^\perp$ the orthogonal complement of $\tilde{B}_{i-1}$ in $ \tilde{B}_i$}
\end{equation*}
relative to the non-degenerate symmetric form $(X,Y)\mapsto \tr_{ \tilde{B}_i/E_i}(XY)$. We also write $1-\alpha:\tilde{A}\rightarrow \tilde{A}$ for the map $x\mapsto x-{}^\alpha x$, whose image is denoted by $\tilde{A}^{-\alpha}$.

\begin{lem}
\label{representatives for X and Y for z}
Each coset  in   $ (s_z^{-1}J_P^-s_z \cap  J_P s_z J_P) / J_P^+$, has a representative $(X,y \varpi_E^{-1}(I+Y'))^+$ with expansions $X=\sum_{j=0}^{d+1}X_j$ and $Y'=\sum_{j=0}^{d+1}Y'_{j}$ such that 
\begin{enumerate}[(i)]
\item $X_0\in \mathfrak{o}_E$ and $X_j\in \tilde{\mathfrak{P}}^{s_{j-1}}_{\Lambda,E_j}$ mod $\tilde{\mathfrak{P}}^{s_{j-1}}_{\Lambda,E_{j-1}} +\tilde{\mathfrak{P}}^{s_{j-1}+1}_{\Lambda,E_j}$, for $j\geq 1$, that lies in $\tilde{\mathfrak{P}}^{s_{j-1}}_{\Lambda,E_j}\cap \tilde{B}_{j-1}^\perp $;

\item $Y_0'\in {\mathfrak{p}}^{}_{E}$ and $Y_j'\in \tilde{\mathfrak{P}}^{s_{j-1}+1}_{\Lambda,E_j}$ for $j\geq 1$, with a decomposition 
$$Y_j'=P_j+Q_j$$ 
where $P_j \in \tilde{\mathfrak{P}}^{s_{j-1}+1}_{\Lambda,E_j} \cap \tilde{B}_{j-1}^\perp$ and $Q_j\in \tilde{\mathfrak{P}}^{r_{j-1}}_{\Lambda,E_j}$, 
satisfying the equations:
\begin{equation}
\label{Y^0_j is in perp}
\sum_{
\begin{smallmatrix}
\max(k,l)=j
\\
k\neq l
\end{smallmatrix}
}X_k{}^\alpha X_l 
 =  (1-\alpha )( y\varpi_E^{-1}P_j).
\end{equation}
and \begin{equation}
\label{solution of Z_j}
X_j{}^\alpha X_j
 =  (1-\alpha )( y\varpi_E^{-1}Q_j),
\end{equation}
(Note that $P_0=0$.)
\end{enumerate}
\end{lem}

\proof We can choose $X_j$, for $j\geq 1$, as stated using the commutative diagram
\begin{equation*}
\label{the commutative diagram
for trace-0-space}
\xymatrixcolsep{3pc}\xymatrix{
\tilde{B}_{j-1}^\perp
\ar@{^{(}->}[r]^{}    
&
\tilde{B}_{j}
\ar[r]^{}  
&
\tilde{B}_{j}/\tilde{B}_{j-1}
\\
\tilde{\mathfrak{P}}^k_{\Lambda,E_j}\cap \tilde{B}_{j-1}^\perp
\ar@{^{(}->}[u]_{}
\ar@{^{(}->}[r]^{ }
&
\tilde{\mathfrak{P}}^k_{\Lambda,E_j}
\ar@{^{(}->}[u]_{}
\ar[r]^{}  
&
\tilde{\mathfrak{P}}^k_{\Lambda,E_j}/\tilde{\mathfrak{P}}^k_{\Lambda,E_{j-1}}.
\ar@{^{(}->}[u]_{}
}
\end{equation*}
We see that the top row is an $\alpha$-equivariant isomorphism 
of $E_j$-spaces, which induces at the bottom row
 an $\alpha$-equivariant isomorphism 
$\tilde{\mathfrak{P}}^k_{\Lambda,E_j}\cap \tilde{B}_{j-1}^\perp \cong \tilde{\mathfrak{P}}^k_{\Lambda,E_j}/\tilde{\mathfrak{P}}^k_{\Lambda,E_{j-1}}
$
of $\mathfrak{o}_{E_j}$-lattices for all $k\in \mathbb{Z}$. For choosing $P_j$, we notice that the left hand side of (\ref{Y^0_j is in perp}) lies in $\tilde{B}_{j-1}^\perp$ by (\ref{trace 0 for X_iX_j}) as well as in the $(-\alpha)$-fixed point subspaces $\tilde{B}_j^{-\alpha}:= \tilde{B}_j\cap \tilde{A}^{-\alpha}$, and the restriction  
$$ 1-\alpha: \tilde{\mathfrak{P}}^{k}_{\Lambda,E_j} \cap \tilde{B}_{j-1}^\perp \rightarrow \tilde{\mathfrak{P}}^{k}_{\Lambda,E_j}\cap \tilde{B}_{j-1}^\perp\cap  \tilde{B}_j^{-\alpha} $$ is surjective. We can also choose $Q_j$ similarly. Now by putting $Y_j'=P_j+Q_j$, we obtain \begin{equation}
\label{Y^0_j is in perp}
\sum_{
\max(k,l)=j
}X_k{}^\alpha X_l 
 =  (1-\alpha )( y\varpi_E^{-1}Y_j'),
\end{equation}
and summing it up for all $j$ yields $X{}^\alpha X = Y-{}^\alpha Y$.
\qed

Here is a very simple consequence which will be frequently used later on: suppose that $T$ is either $X_j$ or $P_j$, for $j\geq 1$, as in Lemma \ref{representatives for X and Y for z}, and $U$ is a product of elements in $B_i$ for $i<j$, then 
\begin{equation}
\label{trace 0 for X_iX_j} 
\mathrm{tr}_{\tilde{B}_j/E_j}(TU)=0.
\end{equation}

\subsubsection{Simplifying the additive parts}
\label{subsection The additive part}

We will simplify the additive parts of $\tilde{\theta}(I+Y')$ and $\theta(I-{}^\alpha XY^{-1}X)$. Since each factor of the additive parts is a value of the character $\psi\circ\mathrm{tr}_{A/F}$, we will show that some of the inputs either lie in $\tilde{\mathfrak{P}}_{\Lambda}$ or have trace 0 by (\ref{trace 0 for X_iX_j}), so that their character values are 1.

\begin{lem}
\label{additive part as quadratic summand 1}
The additive part of $\tilde{\theta}(I+Y')$ is
\begin{equation*}
 \prod_{j=0}^d  \psi\circ \tr_{A/F}\left((c_{j}y_{}^{-1}\varpi_E/2) X_{j+1}{}^\alpha X_{j+1}\right) .
\end{equation*}
\end{lem}

\proof
In our calculation below, we have to switch between the additive expansion of  $Y'\in \tilde{\mathfrak{H}}^1_\Lambda$ given by $Y'=\sum_{i=0}^{d+1}Y'_i$ and the multiplicative expansion according to (\ref{H^1-tilde factorization}). Let's shorthand write 
$$ (I+Y_*)_{0}^{k} = (I+Y_0)\cdots (I+Y_k)
\quad\text{and}\quad
 (I+Y_*)_{k}^{0} = (I+Y_k)\cdots (I+Y_0)$$ for $k=0,\dots,d+1$, so that $Y=y\varpi_E^{-1}(I+Y_*)_{0}^{d+1}$. 
We have also written $Y=y\varpi_E^{-1}(I+Y')$ where $Y'=\sum_{i=0}^{d+1}Y'_i$, such that
\begin{equation}
\label{change Y_j into Y_j'}
Y_i'= (I+Y_*)_{0}^{i-1} Y_i.
\end{equation} The additive part of $\tilde{\theta}(I+Y')$ is therefore equal to
\begin{equation}
\label{additive part tilde}
\prod_{
\begin{smallmatrix}
i,j=0
\\
j<i
\end{smallmatrix}
}^{d+1} \psi_{c_j}(I+Y_i).
\end{equation}
Note that each factor above can be rewritten as 
\begin{equation}
\label{rewrtite Y}
\psi_{c_j}(I+Y_i) = \psi\circ \tr_{A/F}(\frac{c_j}{2}(Y_i+{}^\alpha Y_i)).
\end{equation}
We write $\widehat{Y}_i\in \tilde{\mathfrak{P}}^{s_{i-1}+1}_{\Lambda,E_i}$ such that $(I+Y_i)(I+\widehat{Y}_i)=1$. We fix $j$ and rewrite the factors on the right hand side of (\ref{rewrtite Y}) as 
\begin{equation*}
\label{simplify Y'_j+alpha Y'_j}
\psi\circ \tr_{A/F}(\frac{c_j}{2}((I+\widehat{Y}_*)_{i-1}^{0}Y_i'+{}^\alpha Y_i' (I-{}^\alpha \widehat{Y}_*)_{0}^{i-1})).
\end{equation*}
This expression admits a lot of simplification by writing $Y_i' = P_i+Q_i$ as in Lemma \ref{representatives for X and Y for z}. First of all, we have 
\begin{equation*}
\psi\circ \tr_{A/F}(\frac{c_j}{2}(I+\widehat{Y}_*)_{i-1}^{0} P_i) =1
\end{equation*}
since each summand has trace 0 by (\ref{trace 0 for X_iX_j}), and similarly
\begin{equation*}
\psi\circ \tr_{A/F}(\frac{c_j}{2}{}^\alpha P_i (I-{}^\alpha \widehat{Y}_*)_{0}^{i-1}) ) =1.
\end{equation*}
Moreover, 
\begin{equation*}
\psi\circ \tr_{A/F}(\frac{c_j}{2}(I+\widehat{Y}_*)_{i-1}^{0} Q_i) =\psi\circ \tr_{A/F}(\frac{c_j}{2}Q_i), 
\end{equation*}
since $\frac{c_j}{2}Q_i\in \tilde{\mathfrak{A}}_{\Lambda}$
 and all other summands (involving $\widehat{Y}_*$) lie in $\tilde{\mathfrak{P}}_{\Lambda}$. By the same reason, this term is non-trivial only when $i=j+1$. We have similarly
\begin{equation*}
\psi\circ \tr_{A/F}(\frac{c_j}{2}{}^\alpha Q_i (I-{}^\alpha \widehat{Y}_*)_{0}^{i-1}) ) =
\begin{cases}
\psi\circ \tr_{A/F}(\frac{c_j}{2}{}^\alpha Q_{j+1})&\text{if } i=j+1,
\\
1&\text{otherwise}.
\end{cases}
\end{equation*}
Therefore, (\ref{additive part tilde}) is equal to 
\begin{equation*}
\prod_{j=0}^d
 \psi\circ \tr_{A/F}(\frac{c_j}{2}(Q_{j+1}+{}^\alpha Q_{j+1})) \end{equation*}
and the lemma follows using (\ref{solution of Z_j}).
\qed

\begin{lem}
\label{additive part as quadratic summand 2}
The additive part of $\theta(I-{}^\alpha XY^{-1}X)$ is 
\begin{equation*}
\prod_{k=0}^d\psi\circ\mathrm{tr}_{A/F}\left(-({c_{k}}/{2} ){}^\alpha X_{k+1} 
y^{-1}\varpi_EX_{k+1}\right).
\end{equation*}
\end{lem}

\proof

If we expand $I-{}^\alpha XY^{-1}X = (I+W_*)_0^{d+1}$ using (\ref{H^1 factorization}) for some $W_k\in \tilde{\mathfrak{P}}^{s_{k-1}+1}_{\Lambda,E_k}$, then  
the additive part of 
$\theta(I-{}^\alpha XY^{-1}X)$ can be written as
\begin{equation}
\label{additive part non-tilde}
\prod_{
\begin{smallmatrix}
k,l=0
\\
k<l
\end{smallmatrix}
}^{d+1} \psi_{c_k/2}(I+W_l),
\end{equation}
To express $I-{}^\alpha XY^{-1}X$ additively, we first denote, for every subset $S\subset \{0,\dots,d+1\}$, a shorthand notation
$$(\widehat{Y}_*)_S = \widehat{Y}_{i_1}\cdots \widehat{Y}_{i_{\#S}} \qquad
\text{if }S=\{i_1,\dots,i_{\#S}\}\text{ listed in descending order.}$$
We can expand $I-{}^\alpha XY^{-1}X$  as
\begin{equation*}
\label{1-xyx complete expansion}
I-{}^\alpha XY^{-1}X=I-\sum_{i,j=0}^{d+1}\sum_{S}{}^\alpha X_i 
(\widehat{Y}_*)_S
y^{-1}\varpi_EX_j,
\end{equation*}
so that we can write $I-{}^\alpha XY^{-1}X = I+\sum_{k=0}^{d+1} W_k'$ according to the additive expansion of ${\mathfrak{H}}^1_\Lambda$, where \begin{equation}
\label{introduce W_k'}
W'_k=-\sum_{\max\{i,j,S\}=k} {}^\alpha X_i 
(\widehat{Y}_*)_S
y^{-1}\varpi_E X_j \in \tilde{\mathfrak{P}}^{s_{k-1}+1}_{\Lambda,E_k},
\end{equation}
We hence have
\begin{equation*}
W'_k = (I+W_*)_0^{k-1} W_k.
\end{equation*}
We now fix $k$ and simplify a sub-product from (\ref{additive part non-tilde}):
\begin{equation*}
\prod_{
\begin{smallmatrix}
k<l
\end{smallmatrix}
} \psi_{c_{k}/2}(I+W_l)
\end{equation*}
by first expressing each of its factors as 
\begin{equation*}
\psi_{c_{k}/2}(I+W_{l}) = \psi\circ\tr_{A/F}(\frac{c_{k}}{2}(I+\widehat{W}_{*})_{l-1}^0W_{l}')
\end{equation*}
and using (\ref{introduce W_k'}) to further expand the input for $\psi\circ\mathrm{tr}_{A/F}$ into summands of the form\begin{equation*}
-\frac{c_{k}}{2}(I+\widehat{W}_{*})_{l-1}^0{}^\alpha X_i(\widehat{Y}_*)_Sy^{-1}\varpi_E X_j.
\end{equation*}
The indices for this summand are $(i,j,S)$ such that $\max\{i,j,S\}=l>k$. In the following cases, this $(i,j,S)$-summand has zero trace.
\begin{enumerate}[(i)]
\item If $\max S<l$ and $i\neq j$, then since one of $i$ and $j$ is $l$, the summand has trace 0 by (\ref{trace 0 for X_iX_j}).

\item If $\max S = l$ and exactly one of $i$ and $j$ is also $l$, then the summand lies in $\tilde{\mathfrak{P}}_{\Lambda}$. 

\item If both $i,j<l$, then $\max S = l$. We write the summand as
\begin{equation*}
-\frac{c_{k}}{2}(I+\widehat{W}_{*})_{l-1}^0{}^\alpha X_i\widehat{Y}_l(\widehat{Y}_*)_{S-\{l\}}y^{-1}\varpi_E X_j
\end{equation*}
and further expand $\widehat{Y}_l = \sum_{m=1}^{\infty}(-Y_l)^m$. Any summand involving $(-Y_l)^m$ with $m\geq 2$ lies in $\tilde{\mathfrak{P}}_{\Lambda}$, i.e., we remain to consider summands of the form
\begin{equation*}
\frac{c_{k}}{2}(I+\widehat{W}_{*})_{l-1}^0{}^\alpha X_i{Y}_l(\widehat{Y}_*)_{S-\{l\}}y^{-1}\varpi_E X_j
\end{equation*}
We then change $Y_l$ into $(I+\widehat{Y}_*)_{l-1}^0Y_l'$ using (\ref{change Y_j into Y_j'}) and decompose  $Y_l'=P_l+Q_l$. We then see that any summand involving $P_l$ has trace 0, and any of those involving $Q_l$ lies in $\tilde{\mathfrak{P}}_{\Lambda}$.

\item The remaining case is $i=j=l$. The  summand lies in $\tilde{\mathfrak{P}}_{\Lambda}$ except when $l=k+1$ and $ S =\emptyset$, which is
\begin{equation*}
-\frac{c_{k}}{2}{}^\alpha X_{k+1} y^{-1}\varpi_E X_{k+1} .
\end{equation*}

\end{enumerate}
Therefore, we obtain
\begin{equation*}
\prod_{l>{k}} \psi_{c_{k}/2}(I+W_{l})  =\psi\circ\mathrm{tr}_{A/F}\left(-\frac{c_{k}}{2} {}^\alpha X_{k+1} 
y^{-1}\varpi_EX_{k+1}\right).
\end{equation*}
and the lemma follows by multiplying the above equalities for all $k$ together.
\qed

By Lemmas \ref{additive part as quadratic summand 1} and \ref{additive part as quadratic summand 2}, 
the additive part of    
$\tilde{\theta}(I+Y')\theta(I-{}^\alpha XY^{-1}X)$ is equal to 
\begin{equation}\label{additivepartfinal}
\prod_{j=0}^{d}\psi\circ \tr_{A/F}\left((y^{-1}\varpi_E/2)(c_j X_{j+1}{}- X_{j+1}c_j ){}^\alpha X_{j+1}\right).
\end{equation}
Note that it is independent of the auxiliary $X_0$, as expected from (\ref{Groth-group for X in z}).

\subsubsection{Non-degeneracy of a quadratic form}
\label{subsection Non-degeneracy of a quadratic form}

Define, for $j=0,\dots,d$, a bilinear form on $\mathfrak{W}_{z,j}$ defined in (\ref{coset space W-z}) by
\begin{equation*}
D_j(X,Y) =\tr_{A/F}(y^{-1}\varpi_E
(Xc_j- c_jX) {}^\alpha Y)) \mod \mathfrak{p}_F,\quad \text{for $X,Y\in \mathfrak{W}_{z,j}$}
\end{equation*}
and put
\begin{equation*}
D = D_0\perp\cdots\perp D_d,
\end{equation*}
which defines a bilinear form on $\mathfrak{W}_{z}=\bigoplus_{j=0}^d\mathfrak{W}_{z,j}$  
 such that the decomposition (\ref{coset space W-z}) is orthogonal.

\begin{prop}
\label{non-degenerate quad form}
The quadratic form $D$ is non-degenerate.
\end{prop}
\proof
It suffices to show that each $D_j$ is non-degenerate. Hence we reduce to the situation where 
\begin{itemize}
\item $E/F$ is a tamely ramified extension in $\tilde{A}$, generated by an element $c\in \tilde{A}$ with valuation $-r=-2s-1$,

\item $D(X,Y)=\tr_{A/F}(\varpi_E  a_c(X)Y) \mod \mathfrak{p}_F$ is a bilinear form on $\mathfrak{W}:=\tilde{\mathfrak{P}}^s_{\Lambda}/ \tilde{\mathfrak{P}}^s_{\Lambda,E}+ \tilde{\mathfrak{P}}^{s+1}_{\Lambda}$, where $a_c:X\mapsto Xc-cX$.
\end{itemize}
We want to show that $D$ is non-degenerate, which is equivalent to showing that 
\begin{equation*}
a_c(X)\in \tilde{\mathfrak{P}}^{-s}_{\Lambda}\quad\Rightarrow\quad
X\in \tilde{\mathfrak{P}}^s_{\Lambda,E}+ \tilde{\mathfrak{P}}^{s+1}_{\Lambda}.
\end{equation*}
This is implied by the definition of the critical exponent \cite[(1.4.5)]{BK} 
$$k_0(c,\Lambda)=\max\{k\in \mathbb{Z}: a_c^{-1}\tilde{\mathfrak{P}}^k_{\Lambda}\cap \tilde{\mathfrak{A}}_{\Lambda}\not\subset\tilde{\mathfrak{A}}_{\Lambda,E}+ \tilde{\mathfrak{P}}_{\Lambda}\}$$ 
and the minimality of $c$, i.e., $k_0(c,\Lambda) = v_{\Lambda}(c) = -r$ \cite[(1.4.15)]{BK}.
\qed

\subsubsection{A quadratic Gauss sum}
\label{subsection A quadratic Gauss sum}

  By putting together (\ref{formula for b and c for y }), (\ref{TzXYplus}), Lemmas \ref{propmultcancel}, \ref{additive part as quadratic summand 1}, and \ref{additive part as quadratic summand 2}, we obtain $b_z$, which is
\begin{equation}
\label{full-Gauss-sum}
\begin{split}
\tilde{\rho}(-1)  T_z(s_z)  \sum_{y\in\boldsymbol{\mu}_E}
\tilde{\rho}(y)
\sum_{
X\in \mathfrak{W}_z
}
\prod_{j=0}^{d}\psi\circ \tr_{A/F}\left((y^{-1}\varpi_E/2)(c_j X_{j+1}- X_{j+1}c_j ){}^\alpha X_{j+1}\right)
\end{split}
\end{equation}
For a fixed $y\in\boldsymbol{\mu}_E$, the above inner sum over $X\in \mathfrak{W}_z$
 is a quadratic Gauss sum, defined on the $\mathbf{k}_F$-space $\mathfrak{W}_z$ equipped with the quadratic form in Proposition \ref{non-degenerate quad form} which is non-degenerate. Note that, since $\Lambda$ is an $\mathfrak{o}_E$-lattice chain, $\mathfrak{W}_z$ is a $\mathbf{k}_E$-space as well. From basic properties of quadratic Gauss sums, we can express the inner sum as
 \begin{equation*}
\left(\frac{2y^{-1}}{\boldsymbol{\mu}_E}\right)^{\dim_{\mathbf{k}_E}\mathfrak{W}_z}
\sum_{X\in \mathfrak{W}_z}\prod_{j=0}^{d}
{\psi}\circ \tr_{A/F}(\varpi_E(c_j X_{j+1}- X_{j+1}c_j ){}^\alpha X_{j+1}).
 \end{equation*}
The sum
\begin{equation}
\label{the-quad-Gauss-sum}
\sum_{X\in \mathfrak{W}_z}\prod_{j=0}^{d}
{\psi}\circ \tr_{A/F}(\varpi_E(c_j X_{j+1}- X_{j+1}c_j ){}^\alpha X_{j+1}),
\end{equation} 
 is equal to a 4-th root of unity, denoted by $\epsilon_z^P(\varpi_E,\mathbf{s},\psi)$, times the positive number $(\#\mathfrak{W}_z)^{1/2}$. Now note that $\dim_{\mathbf{k}_E}\mathfrak{W}_z=fe-1$ where $e=e(E/F)$, $f=f(E/F)$ and, in our situation when $q$ and $e$ are odd, since $$\epsilon_z^P(\varpi_E,\mathbf{s},\psi)^2=
\left(\frac{-1}{\boldsymbol\mu_E}\right)^{\dim_{\mathbf{k}_E}\mathfrak{W}_z}
=
\left(\frac{-1}{\boldsymbol\mu_E}\right)^{f-1}
=(-1)^{(q^f-1)(f-1)/2}=1,$$ 
the normalized sum  $\epsilon_z^P(\varpi_E,\mathbf{s},\psi)$ is actually a sign.

Therefore, 
\begin{equation*}
b_z=0\quad \text{ when }\quad\tilde{\rho}|_{\boldsymbol{\mu}_E}=\left(\frac{\cdot}{\boldsymbol\mu_E}\right)^f.
\end{equation*}
and 
\begin{equation*}
b_z=\tilde{\rho}(-2)\epsilon_z^P(\varpi_E,\mathbf{s},\psi)
T_z(s_z)\#(\mathfrak{W}_z)^{1/2}(q_E-1) \quad\text{ when }\quad \tilde{\rho}|_{\boldsymbol{\mu}_E}=\left(\frac{\cdot}{\boldsymbol\mu_E}\right)^{f-1}.
\end{equation*}
By noting that $\#(\mathfrak{W}_z) = c_z / q_E$, we finish the proof of Theorem \ref{eigenvalues of Ty and Tz}.(ii).

Finally, we show that our main results are independent of the choices of the additive character $\psi$ and the uniformizer $\varpi_E$,
as they should be.
\begin{prop}
\begin{enumerate}[(i)]
    \item The sign $\epsilon_z^P(\varpi_E,\mathbf{s},\psi)$ is independent of the additive character $\psi$, hence denoted by $\epsilon_z^P(\varpi_E,\mathbf{s})$. 
    
    \item The relation $\tilde{\rho}(\varpi_E)=  \epsilon_z^P(\varpi_E,\mathbf{s})\rho(-1)$ in Theorem \ref{main theorem} is independent of the choice of $\varpi_E$. 
\end{enumerate}
\end{prop}

\proof
For (i), show that the sign $\epsilon_z^P(\varpi_E,\mathbf{s},\psi)$ is independent of the chosen additive character $\psi$ of $F$. Indeed if we replace $\psi$ by $\psi_a:F\rightarrow \mathbb{C}^\times,\, x\mapsto \psi(ax)$ for some $a\in F^\times$, then the Gauss sum is multiplied by $\left(\frac{a}{\boldsymbol\mu_E}\right)^{f-1}$. Note that the character $\left(\frac{\cdot}{\boldsymbol\mu_E}\right)^{f-1}$ is non-trivial only when $f$ is even, so that $[\boldsymbol\mu_E:\boldsymbol\mu_F] = (q^f-1)/(q-1)$ is also even and $\left(\frac{\cdot}{\boldsymbol\mu_E}\right)^{f-1}$ is always trivial on $\boldsymbol\mu_F$. To prove (ii), suppose that $\varpi_E$ is replaced by $u\varpi_E$ for some $u\in U_E$. It is easy to see that $
\epsilon_z^P(u \varpi_E,\mathbf{s}) = \tilde{\rho}(u) \epsilon_z^P(\varpi_E,\mathbf{s})
$
from the expression (\ref{full-Gauss-sum}). 
\qed

\addcontentsline{toc}{section}{References}

\def\cprime{$'$} \def\cftil#1{\ifmmode\setbox7\hbox{$\accent"5E#1$}\else
  \setbox7\hbox{\accent"5E#1}\penalty 10000\relax\fi\raise 1\ht7
  \hbox{\lower1.15ex\hbox to 1\wd7{\hss\accent"7E\hss}}\penalty 10000
  \hskip-1\wd7\penalty 10000\box7}

\end{document}